\documentclass[preprint,12pt]{elsarticle}

\usepackage{hyperref,amsmath,amsfonts,mathrsfs}
\journal{}

\begin{document}

\begin{frontmatter}



\title{Spectra and energy of bipartite signed digraphs}


\author{Mushtaq A. Bhat and S. Pirzada}

\address{Department of Mathematics, University of Kashmir, Srinagar, India\\mushtaqab1125@gmail.com; pirzadasd@kashmiruniversity.ac.in}

\begin{abstract}
The set of distinct eigenvalues of a signed digraph $S$ together with their multiplicities is called its spectrum. The energy of a signed digraph $S$ with eigenvalues $z_1,z_2,\cdots,z_n$ is defined as $E(S)=\sum_{j=1}^{n}|\Re z_j|$, where $\Re z_j $ denotes real part of complex number $z_j$. In this paper, we show that the characteristic polynomial of a bipartite signed digraph of order $n$ with each cycle of length $\equiv 0\pmod 4$ negative and each cycle of length  $\equiv 2\pmod 4$ positive is of the form \\
$$\phi_S(z)=z^n+\sum\limits_{j=1}^{\lfloor{\frac{n}{2}}\rfloor}(-1)^j c_{2j}(S)z^{n-2j},$$\\
where $c_{2j}(S)$ are nonnegative integers. We define a quasi-order relation in this case and show energy is increasing. It is shown that the characteristic polynomial of a bipartite signed digraph of order $n$ with each cycle negative has the form \\
$$\phi_S(z)=z^n+\sum\limits_{j=1}^{\lfloor{\frac{n}{2}}\rfloor}c_{2j}(S)z^{n-2j},$$
where $c_{2j}(S)$ are nonnegative integers. We study integral, real, Gaussian signed digraphs and quasi-cospectral digraphs and show for each positive integer $n\ge 4$ there exists a family of $n$ cospectral, non symmetric, strongly connected, integral, real, Gaussian signed digraphs (non cycle balanced) and quasi-cospectral digraphs of order $4^n$.  We obtain a new family of pairs of equienergetic strongly connected signed digraphs and answer to open problem $(2)$ posed in Pirzada and Mushtaq, Energy of signed digraphs, Discrete Applied Mathematics 169 (2014) 195-205.
\end{abstract}


\end{frontmatter}

% \linenumbers

\section{Introduction}\label{sec1}

A signed digraph (or briefly sidigraph) is defined to be a pair $S=(D,\sigma)$ where $D=(V,\mathscr{A})$ is the underlying digraph and $\sigma: \mathscr{A}\rightarrow \{-1, 1\}$ is the signing function. The sets of positive and negative arcs of $S$ are respectively denoted by $\mathscr{A}^+(S)$ and $\mathscr{A}^-(S)$. Thus $\mathscr{A}(S)=\mathscr{A}^+(S)\cup \mathscr{A}^- (S)$. A sidigraph is said to be homogeneous if all of its arcs have either positive sign or negative sign, and heterogeneous otherwise. Throughout this paper, bold arcs will denote positive arcs and dotted arcs will denote negative arcs. \\
\indent Two vertices are adjacent if they are connected by an arc. If there is an arc from a vertex $u$ to the vertex $v$, we indicate this by $(u,v)$. A path of length $n-1$ $(n\geq2)$, denoted by $P_n$, is a sidigraph  on $n$ vertices $v_1, v_2,\cdots,v_n$ with $n-1$  signed arcs $(v_i,v_{i+1})$, $i=1,2,\cdots,n-1$. A cycle of length $n$ is a sidigraph having vertices $v_1, v_2,\cdots,v_n$ and signed arcs $(v_i, v_{i+1})$, $i=1,2,\cdots,n-1$ and $(v_n,v_1)$. A sidigraph $S$ is said to be strongly connected if its underlying digraph $S^u$ is strongly connected. The sign of a sidigraph is defined as the product of signs of its arcs. A sidigraph is said to be positive (negative) if its sign is positive (negative) i.e., it contains an even (odd) number of negative arcs. A sidigraph is said to be all-positive (respectively, all-negative) if all its arcs are positive (negative). A sidigraph is said to be cycle balanced if each of its cycles is positive, otherwise non cycle balanced. The negative of a sidigraph $S$ denoted by $-S$ is the sidigraph obtained by negating sign of each arc of $S$. Throughout we call cycle balanced cycle a positive cycle and non cycle balanced cycle a negative cycle and respectively denote them by $C_n$ and ${\bf C}_n$, where $n$ is number of vertices.\\
\indent The adjacency matrix of a sidigraph $S$ with vertex set $\{v_1,v_2,\cdots,v_n\}$ is the $n\times n$ matrix $A(S)=(a_{ij})$, where\\
$$a_{ij}=\left\{\begin{array}{lr}\sigma(v_i,v_j), &\mbox{if there is an arc from $v_i$ to $v_j,$}\\
0, &\mbox {otherwise.}
\end{array}\right.$$ \\
\indent The characteristic polynomial $|zI-A(S)|$ of the adjacency matrix $A(S)$  of a sidigraph $S$ is called the characteristic polynomial of $S$ and is denoted by $\phi_S(z)$. The eigenvalues of $A(S)$ are called the eigenvalues of $S$. The set of distinct eigenvalues of $S$ together with their multiplicities is called the spectrum of $S$. If $S$ is a sidigraph of order $n$ with distinct eigenvalues $z_1,z_2,\cdots,z_k$ and if their respective multiplicities are $m_1,m_2,\cdots,m_k$, we write the spectrum of $S$ as $spec(S)=\{z^{(m_1)}_1,z^{(m_2)}_2,\cdots,z^{(m_k)}_k\}$.\\
\indent A sidigraph is symmetric if $(u,v)\in \mathscr{A}^{+}(S)$ ~(or~$\mathscr{A}^-(S)$) then $(v,u)\in \mathscr{A}^{+} (S)$~(or~~$\mathscr{A}^-(S)$), where $u,v\in V(S)$. A one to one correspondence between sigraphs and symmetric sidigraphs is given by $\Sigma\rightsquigarrow \overleftrightarrow{\Sigma}$, where $\overleftrightarrow{\Sigma}$ has the same vertex set as that of sigraph $\Sigma$, and each signed edge $(u,v)$ is replaced by a pair of symmetric arcs $(u,v)$ and $(v,u)$ both with same sign as that of edge $(u,v)$. Under this correspondence a sigraph can be identified with a symmetric sidigraph. A sidigraph is said to be skew symmetric if its adjacency matrix is skew symmetric. A linear subsidigraph of a sidigraph $S$ is a subsidigraph with indegree and out degree of each vertex equal to one i.e., each component is a cycle.\\

\indent The following is the coefficient theorem for sidigraphs \cite{a}.\\

\noindent {\bf Theorem 1.1.} If $S$ is a sidigraph with characteristic polynomial
$$\phi_S(z)=z^n+b_1(S) z^{n-1} +\cdots+b_{n-1}(S)z+b_n(S)$$\\
then $$ b_j(S)=\sum\limits_{L\in \pounds_j}(-1)^{p(L)}\prod\limits_{Z\in c(L)} s(Z),$$
for all $j=1,2,\cdots,n$, where $\pounds_j$ is the set of all linear subsidigraphs  $L$ of $S$ of order $j$, $p(L)$ denotes number of components of $L$ and $c(L)$ denotes the set of all cycles of $L$ and $s(Z)$ the sign of cycle $Z$.\\

\indent The spectral criterion for cycle balance of sidigraphs (sigraphs) given by Acharya \cite{a} is as follows.\\

\noindent {\bf Theorem 1.2.} A sidigraph (sigraph) $S$ is cycle balanced (balanced) if and only if it is cospectral with the underlying unsigned digraph (graph).\\

\indent The energy of a graph is defined as the sum of the absolute values of graph eigenvalues. This  concept was given by Gutman \cite{g}. Later the concept of energy was extended to digraphs by Pe\~{n}a and Rada and they defined the energy of a digraph to be the sum of the absolute values of the real parts of digraph eigenvalues \cite{pr}. The energy of a sigraph $\Sigma$ was defined by Germina, Hameed and Zaslavsky \cite{ghz} as $E(\Sigma)=\sum_{j=1}^{n}|\lambda_j|,$ where $\lambda_1,\lambda_2,\cdots,\lambda_n$ are eigenvalues of the sigraph $\Sigma$. Recently, this concept was generalized to sidigraphs \cite{pm} as $E(S)=\sum_{j=1}^{n}|\Re z_j|,$ where $z_1,z_2,\cdots,z_n$ are the eigenvalues of $S$ and $\Re z_j$ denotes real part of $z_j$. This definition was motivated by Coulson's integral formula for sidigraphs [Theorem $3.3$, \cite{pm}]\\
$$ E(S)=\sum\limits_{j=1}^{n}|\Re z_j|=\frac{1}{\pi} \int\limits_{-\infty}^{\infty}(n-\frac{\iota x \phi'_{S}(\iota x)}{\phi_{S}(\iota x)})dx$$ \\
where $\iota=\sqrt{-1}$ and  $\int\limits_{-\infty}^{\infty} F(x) dx$ denotes the principle value of the respective integral. We note that in \cite{rb} Brualdi calls this type of energy as the low energy.\\

For more information about spectra and energy of graphs, sigraphs, digraphs and sidigraphs see \cite{rb,cds,g,g1,pm,pr,r,r1}.

\indent Esser and Harary in \cite{eh2} showed that the spectrum of a strongly connected digraph $D$ remains invariant under the multiplication by $-1$ if and only if it is bipartite. We show there are non bipartite strongly connected sidigraphs with this property. As in bipartite digraphs, in general, the even coefficients of non cycle balanced bipartite sidigraphs does not alternate in sign. For example the characteristic polynomial of non cycle balanced bipartite sidigraph $S$ in Fig. $3$ is $\phi_S(z)=z^4+z^2$. Clearly, even coefficients do not alternate in sign. Now, consider the non cycle balanced bipartite sidigraph $S_1$ in Fig. $2$, the characteristic polynomial is $\phi_{S_1}(z)=z^6-z^4+2z^2$. In this case, even coefficients alternate in sign. In \cite{rgc}, the authors considered the bipartite digraphs with characteristic polynomial of the form\\
\begin{equation}
\phi_D(z)=z^n+\sum\limits_{j=1}^{\lfloor{\frac{n}{2}}\rfloor}(-1)^j c_{2j}(D)z^{n-2j}
\end{equation}
where $c_{2j}(D)$ are nonnegative integers for every $j=1,2,\cdots,{\lfloor{\frac{n}{2}}\rfloor}$ and studied a large family of bipartite digraphs $\Delta_n$ consisting of digraphs with $n$ vertices where each cycle has length $\equiv 2\pmod 4$ with characteristic polynomial of the form $(1)$. Because of this alternating nature of even coefficients it was possible to compare energies of digraphs in $\Delta_n$ by means of quasi-order relation. It is natural to consider the same problem for sidigraphs. We show bipartite sidigraphs on $n$ vertices with each cycle of length $\equiv 0\pmod 4$ negative (i.e., containing odd number of negative arcs) and each cycle of length $\equiv 2\pmod 4$ positive (i.e., containing an even number of negative arcs) has characteristic polynomial of the form $(1)$. We denote this class of sidigraphs by $\Delta^1_n$. We derive an integral expression for the energy and define a quasi-order relation to compare the energies of sidigraphs in this case. We also study another class of bipartite sidigraphs on $n$ vertices with all cycles negative (i.e., each cycle has odd number of negative arcs) and show a sidigraph in this class has characteristic polynomial of the form \\
\begin{equation}
\phi_S(z)=z^n+\sum\limits_{j=1}^{\lfloor{\frac{n}{2}}\rfloor}c_{2j}(S)z^{n-2j}
\end{equation}
where $c_{2j}(S)$ are nonnegative integers for every $j=1,2,\cdots,{\lfloor{\frac{n}{2}}\rfloor}$. We denote this class of sidigraphs by $\Delta^2_n$.

\indent Two sidigraphs of same order are said to be cospectral (or isospectral) if they have the same spectrum. In \cite{eh1}, the authors studied digraphs with integral, real and Gaussian spectra. We study sidigraphs with integral, real and Gaussian spectra and we show for each positive integer $n\ge 4$ there exists a collection of $n$ non cycle balanced, non symmetric, strongly connected, integral, real and Gaussian cospectral sidigraphs of order $4^n$.\\
\indent Two noncospectral sidigraphs of same order are said to be equienergetic if they have the same energy. Equienergetic sidigraphs were obtained in \cite{pm} and the authors raised the following open problem.\\

\noindent {\bf Problem 1.3.} Find an infinite family of pairs of noncospectral equienergetic sidigraphs on $n\ge 4$ vertices with both constituents non cycle balanced.\\

In this paper, we obtain a new family of strongly connected equienergetic sidigraphs and give the answer to Problem $1.3$.\\

\section{Spectra of sidigraphs}
Recall a sidigraph $S$ is bipartite if its underlying digraph is bipartite. The following result by Esser and Harary \cite{eh2} characterizes strongly connected bipartite digraphs in terms of spectra.\\

\noindent {\bf Theorem 2.1.} A strongly connected digraph $D$ is bipartite if and only if its spectrum is invariant under multiplication by $-1$.\\

\indent Let $S$ be a bipartite sidigraph, then as in [Theorem $3.3$, \cite{cds}], the characteristic polynomial of $S$ is given by $\phi_S(z)=z^{\delta}\psi(z^2)$, where $\delta$ is a nonnegative integer and $\psi(z^2)$ is a polynomial in $z^2$. Therefore the spectrum of a bipartite sidigraph remains invariant under multiplication by $-1$.\\

%TeXCAD Picture [last.pic]. Options:
%\grade{\on}
%\emlines{\off}
%\epic{\off}
%\beziermacro{\on}
%\reduce{\on}
%\snapping{\off}
%\quality{8.00}
%\graddiff{0.01}
%\snapasp{1}
%\zoom{4.0000}
\unitlength 1mm % = 2.85pt
\linethickness{0.4pt}
\ifx\plotpoint\undefined\newsavebox{\plotpoint}\fi % GNUPLOT compatibility
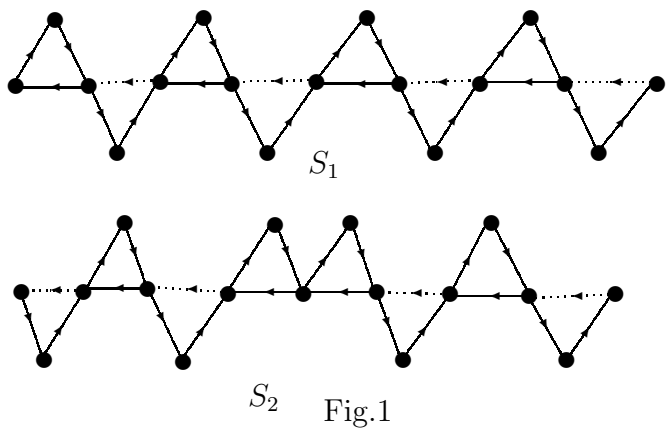
\begin{picture}(96.5,68.87)(0,0)
\put(10.5,58.75){\circle*{2}}
\put(20.25,58.75){\circle*{2}}
\put(29.75,59.25){\circle*{2}}
\put(39.25,59){\circle*{2}}
\put(50.5,59.25){\circle*{2}}
\put(61.5,59){\circle*{2}}
\put(72.25,59){\circle*{2.}}
\put(83.5,59){\circle*{2}}
\put(95.75,59){\circle*{2}}
\put(15.75,67.5){\circle*{2}}
\put(24,49.75){\circle*{2}}
\put(35.5,67.75){\circle*{2}}
\put(44,49.75){\circle*{2}}
\put(57.25,67.75){\circle*{2}}
\put(66.25,49.75){\circle*{2}}
\put(79,67.75){\circle*{2}}
\put(88,49.75){\circle*{2}}
\put(11.25,31.25){\circle*{2}}
\put(19.5,31.25){\circle*{2}}
\put(38.75,31){\circle*{2}}
\put(48.75,31){\circle*{2}}
\put(68.25,31){\circle*{2}}
\put(78.75,30.75){\circle*{2}}
\put(90.25,31){\circle*{2}}
\put(14.25,22.25){\circle*{2}}
\put(25,40.5){\circle*{2}}
\put(32.75,22){\circle*{2}}
\put(45,40.25){\circle*{2}}
\put(55,40.5){\circle*{2}}
\put(62,22.25){\circle*{2}}
\put(73.75,40.5){\circle*{2}}
\put(83.75,22.25){\circle*{2}}
%\vector[middle](10.25,59)(15.75,68.25)
\put(13,63.63){\vector(2,3){.07}}\multiput(10.25,59)(.03353659,.05640244){164}{\line(0,1){.05640244}}
%\end
%\vector[middle](15.75,67.75)(20.25,58.5)
\put(18,63.13){\vector(1,-2){.07}}\multiput(15.75,67.75)(.03358209,-.06902985){134}{\line(0,-1){.06902985}}
%\end
%\vector[middle](20,58.5)(10,58.5)
\put(15,58.5){\vector(-1,0){.07}}\put(20,58.5){\line(-1,0){10}}
%\end
%\vector[middle](20.25,59)(24.25,49.75)
\put(22.25,54.38){\vector(1,-2){.07}}\multiput(20.25,59)(.03361345,-.07773109){119}{\line(0,-1){.07773109}}
%\end
%\vector[middle](24,50.25)(29.5,59.5)
\put(26.75,54.88){\vector(2,3){.07}}\multiput(24,50.25)(.03353659,.05640244){164}{\line(0,1){.05640244}}
%\end
%\vector[middle](29.25,59.25)(35.25,68.25)
\put(32.25,63.75){\vector(2,3){.07}}\multiput(29.25,59.25)(.03370787,.0505618){178}{\line(0,1){.0505618}}
%\end
%\vector[middle](35.25,68.25)(39.25,59)
\put(37.25,63.63){\vector(1,-2){.07}}\multiput(35.25,68.25)(.03361345,-.07773109){119}{\line(0,-1){.07773109}}
%\end
%\vector[middle](39.25,59)(43.75,50)
\put(41.5,54.5){\vector(1,-2){.07}}\multiput(39.25,59)(.03358209,-.06716418){134}{\line(0,-1){.06716418}}
%\end
%\vector[middle](44,49.75)(50.75,59.5)
\put(47.38,54.63){\vector(2,3){.07}}\multiput(44,49.75)(.03358209,.04850746){201}{\line(0,1){.04850746}}
%\end
%\vector[middle](50.5,59.5)(57,68)
\put(53.75,63.75){\vector(3,4){.07}}\multiput(50.5,59.5)(.03367876,.04404145){193}{\line(0,1){.04404145}}
%\end
%\vector[middle](57,68)(61.5,59.25)
\put(59.25,63.63){\vector(1,-2){.07}}\multiput(57,68)(.03358209,-.06529851){134}{\line(0,-1){.06529851}}
%\end
%\vector[middle](61.5,58.75)(65.75,50)
\put(63.63,54.38){\vector(1,-2){.07}}\multiput(61.5,58.75)(.03373016,-.06944444){126}{\line(0,-1){.06944444}}
%\end
%\vector[middle](66,50)(72.25,59.5)
\put(69.13,54.75){\vector(2,3){.07}}\multiput(66,50)(.03360215,.05107527){186}{\line(0,1){.05107527}}
%\end
%\vector[middle](72.25,59.5)(79,68.25)
\put(75.63,63.88){\vector(3,4){.07}}\multiput(72.25,59.5)(.03358209,.04353234){201}{\line(0,1){.04353234}}
%\end
%\vector[middle](79,68.25)(83,59.5)
\put(81,63.88){\vector(1,-2){.07}}\multiput(79,68.25)(.03361345,-.07352941){119}{\line(0,-1){.07352941}}
%\end
%\vector[middle](83.5,59.25)(87.75,50.25)
\put(85.63,54.75){\vector(1,-2){.07}}\multiput(83.5,59.25)(.03373016,-.07142857){126}{\line(0,-1){.07142857}}
%\end
%\vector[middle](87.75,49.75)(95.25,59.25)
\put(91.5,54.5){\vector(3,4){.07}}\multiput(87.75,49.75)(.03363229,.0426009){223}{\line(0,1){.0426009}}
%\end
%\vector[middle](83.5,59.25)(72.25,59.25)
\put(77.88,59.25){\vector(-1,0){.07}}\put(83.5,59.25){\line(-1,0){11.25}}
%\end
%\vector[middle](61.5,59)(49.75,59)
\put(55.63,59){\vector(-1,0){.07}}\put(61.5,59){\line(-1,0){11.75}}
%\end
%\vector[middle](39.25,59)(29.5,59)
\put(34.38,59){\vector(-1,0){.07}}\put(39.25,59){\line(-1,0){9.75}}
%\end
%\vector[middle]{dot}(95.75,59.25)(83.5,59.25)
\put(89.63,59.25){\vector(-1,0){.07}}\multiput(95.68,59.18)(-.9423,0){14}{{\rule{.4pt}{.4pt}}}
%\end
%\vector[middle]{dot}(71.75,59.5)(61.25,59)
\put(66.5,59.25){\vector(-1,0){.07}}\multiput(71.68,59.43)(-.9545,-.0455){12}{{\rule{.4pt}{.4pt}}}
%\end
%\vector[middle]{dot}(50.25,59.5)(39,59.25)
\put(44.63,59.38){\vector(-1,0){.07}}\multiput(50.18,59.43)(-.9375,-.0208){13}{{\rule{.4pt}{.4pt}}}
%\end
%\vector[middle]{dot}(29.5,59.5)(20,59)
\put(24.75,59.25){\vector(-1,0){.07}}\multiput(29.43,59.43)(-.95,-.05){11}{{\rule{.4pt}{.4pt}}}
%\end
\put(51.5,48.25){\makebox(0,0)[cc]{$S_1$}}
\put(43.5,17.25){\makebox(0,0)[cc]{$S_2$}}
\put(56,15.25){\makebox(0,0)[cc]{Fig.$1$}}
%\vector[middle](11,31.75)(14,22.75)
\put(12.5,27.25){\vector(1,-3){.07}}\multiput(11,31.75)(.0337079,-.1011236){89}{\line(0,-1){.1011236}}
%\end
%\vector[middle](14,22.5)(19.5,31.25)
\put(16.75,26.88){\vector(2,3){.07}}\multiput(14,22.5)(.03353659,.05335366){164}{\line(0,1){.05335366}}
%\end
%\vector[middle](19.75,31.5)(25,40.75)
\put(22.38,36.13){\vector(1,2){.07}}\multiput(19.75,31.5)(.03365385,.05929487){156}{\line(0,1){.05929487}}
%\end
%\vector[middle](24.75,41)(28.25,31)
\put(26.5,36){\vector(1,-3){.07}}\multiput(24.75,41)(.03365385,-.09615385){104}{\line(0,-1){.09615385}}
%\end
%\vector[middle](32.5,22.25)(38.75,31.25)
\put(35.63,26.75){\vector(2,3){.07}}\multiput(32.5,22.25)(.03360215,.0483871){186}{\line(0,1){.0483871}}
%\end
%\vector[middle](38.75,31.5)(44.75,40.75)
\put(41.75,36.13){\vector(2,3){.07}}\multiput(38.75,31.5)(.03370787,.05196629){178}{\line(0,1){.05196629}}
%\end
%\vector[middle](48.5,31.25)(55.25,40.75)
\put(51.88,36){\vector(3,4){.07}}\multiput(48.5,31.25)(.03358209,.04726368){201}{\line(0,1){.04726368}}
%\end
%\vector[middle](55.25,40.25)(58.25,31.25)
\put(56.75,35.75){\vector(1,-3){.07}}\multiput(55.25,40.25)(.0337079,-.1011236){89}{\line(0,-1){.1011236}}
%\end
%\vector[middle](58.5,31.25)(61.75,22.25)
\put(60.13,26.75){\vector(1,-3){.07}}\multiput(58.5,31.25)(.03350515,-.09278351){97}{\line(0,-1){.09278351}}
%\end
%\vector[middle](62,22.5)(68.25,31)
\put(65.13,26.75){\vector(3,4){.07}}\multiput(62,22.5)(.03360215,.04569892){186}{\line(0,1){.04569892}}
%\end
%\vector[middle](68.25,31)(73.5,41.25)
\put(70.88,36.13){\vector(1,2){.07}}\multiput(68.25,31)(.03365385,.06570513){156}{\line(0,1){.06570513}}
%\end
%\vector[middle](73.75,40.75)(78.75,31)
\put(76.25,35.88){\vector(1,-2){.07}}\multiput(73.75,40.75)(.03355705,-.06543624){149}{\line(0,-1){.06543624}}
%\end
%\vector[middle](78.75,30.5)(83.5,22.25)
\put(81.13,26.38){\vector(1,-2){.07}}\multiput(78.75,30.5)(.03368794,-.05851064){141}{\line(0,-1){.05851064}}
%\end
%\vector[middle](83.75,22.25)(89.75,31)
\put(86.75,26.63){\vector(2,3){.07}}\multiput(83.75,22.25)(.03370787,.0491573){178}{\line(0,1){.0491573}}
%\end
%\vector[middle](78.75,30.75)(68.25,30.75)
\put(73.5,30.75){\vector(-1,0){.07}}\put(78.75,30.75){\line(-1,0){10.5}}
%\end
%\vector[middle]{dot}(90.25,31)(78.75,30.5)
\put(84.5,30.75){\vector(-1,0){.07}}\multiput(90.18,30.93)(-.9583,-.0417){13}{{\rule{.4pt}{.4pt}}}
%\end
%\vector[middle]{dot}(68,31)(58.25,31.25)
\put(63.13,31.13){\vector(-1,0){.07}}\multiput(67.93,30.93)(-.975,.025){11}{{\rule{.4pt}{.4pt}}}
%\end
%\vector[middle]{dot}(19.25,31.5)(11,31.5)
\put(15.13,31.5){\vector(-1,0){.07}}\multiput(19.18,31.43)(-.9167,0){10}{{\rule{.4pt}{.4pt}}}
%\end
\put(28,31.75){\circle*{2}}
%\vector[middle](58.25,31.25)(48.5,31.25)
\put(53.38,31.25){\vector(-1,0){.07}}\put(58.25,31.25){\line(-1,0){9.75}}
%\end
%\vector[middle](48.5,31.25)(38.5,31.25)
\put(43.5,31.25){\vector(-1,0){.07}}\put(48.5,31.25){\line(-1,0){10}}
%\end
%\vector[middle](45,40.75)(48.5,31.5)
\put(46.75,36.13){\vector(1,-3){.07}}\multiput(45,40.75)(.03365385,-.08894231){104}{\line(0,-1){.08894231}}
%\end
\put(58.5,31.25){\circle*{2}}
%\vector[middle](27.75,31.75)(19.25,31.75)
\put(23.5,31.75){\vector(-1,0){.07}}\put(27.75,31.75){\line(-1,0){8.5}}
%\end
%\vector[middle](28,32)(32.5,22.5)
\put(30.25,27.25){\vector(1,-2){.07}}\multiput(28,32)(.03358209,-.07089552){134}{\line(0,-1){.07089552}}
%\end
%\vector[middle]{dot}(38.75,31.25)(27.75,32)
\put(33.25,31.63){\vector(-1,0){.07}}\multiput(38.68,31.18)(-.9167,.0625){13}{{\rule{.4pt}{.4pt}}}
%\end
\end{picture}

\noindent {\bf Remark 2.2.} Unlike in digraphs, the converse of Theorem $2.1$ is not true for sidigraphs. For example sidigraphs $S_1$ and $S_2$ in Fig. 1 are two strongly connected non bipartite sidigraphs of order $17$. It is easy to check that $\phi_{S_1}(z)=\phi_{-S_1}(z)=z^{17}+3z^{11}+z^5$ and $\phi_{S_2}(z)=\phi_{-S_2}(z)=z^{17}+z^{11}+z^5$.\\
\indent We recall the definitions of type $a$, type $b$, type $c$ and type $d$ linear sidigraphs from \cite{pm}. From Theorem $1.1$,  $b_j(S)=\sum\limits_{L\in \pounds_j } (-1)^{p(L)}s(L), j=1,2,\cdots,n$, where $s(L)=\prod\limits_{Z\in c(L)} s(Z)$. Clearly, this sum contains positive and negative ones.\\
$+1$ occurs if, and only if\\ $(a)$ Number of components of $L\in \pounds_j$ is odd  and $s(L)<0$. We call such linear sidigraphs as type $a$ linear sidigraphs.\\
$(b)$ Number of components of $L\in \pounds_j$ is even and $s(L)> 0$. We call such linear sidigraphs as type $b$.\\
$-1$ will occur if, and only if\\ $(c)$ Number of components of $L\in \pounds_j$ is odd and $s(L)> 0$. We call such linear sidigraphs as type $c$.\\ $(d)$ Number of components of $L\in \pounds_j$ is even and $s(L)< 0$. We call such linear sidigraphs as type $d$.\\

We have the following result.\\

\noindent {\bf Theorem 2.3.} Let $S$ be a sidigraph of order $n$. Then the following statements are equivalent.\\
(I) Spectrum of $S$ remains invariant under multiplication by $-1$.\\
(II) $S$ and $-S$ are cospectral, where $-S$ is the negative of $S$.\\
(III) In $S$, for each odd $j$, the number of linear subsidigraphs of order $j$ of type $a$ or type $b$ or both types is equal to the number of linear subsidigraphs of order $j$ of type $c$ or type $d$ or both types.\\
\noindent {\bf Proof.} $(I)\Longrightarrow(II)$ Follows by the fact that $spec(-A)=-spec(A)$ for any square matrix $A$.\\
$(II)\Longrightarrow (III)$ Assume $S$ and $-S$ are cospectral. Then $\phi_{S}(z)=\phi_{-S}(z)=\pm\phi_S(-z)$. The sign is positive or negative according as  the order $n$ of sidigraph is respectively even or odd. This clearly indicates that the coefficient $b_j(S)=0$  for each odd j and therefore $(III)$ follows.\\
$(III)\Longrightarrow (I)$ is immediate.\qed\\

\noindent {\bf Lemma 2.4.} If $S$ is a bipartite sidigraph, then for all $j=1,2,\cdots$\\
$(I)$ $\pounds_{2j-1}=\emptyset$.\\
$(II)$ Every element of $\pounds_{4j}$ has an even number of cyclic components of length $\equiv 2\pmod 4$. The number of components of length $\equiv 0\pmod 4$ is either even or odd.\\
$(III)$ Every element of $\pounds_{4j+2}$ has an odd number of cyclic components of length $\equiv 2\pmod 4$. The number of components of length $\equiv 0\pmod 4$ is either even or odd.\\
\noindent {\bf Proof.} $(I)$. Since $S$ is bipartite, therefore $\pounds_{2j-1}=\emptyset$ for all $j=1,2,\cdots$.\\
$(II)$. Assume $L\in\pounds_{4j}$ has $p$ components of length $4l_r+2$, for $r=1,2,\cdots,p$ and $q$ components of length $4m_r$ for $r=1,2,\cdots,q$. Then
$$4j=\sum\limits_{r=1}^{p}(4l_r+2)+\sum\limits_{r=1}^{q}(4m_r)$$\\
which gives $p=2j-2\sum\limits_{r=1}^{p}(l_r)-2\sum\limits_{r=1}^{q}(m_r)$. This shows $p$ is even irrespective of whether $q$ is even or odd.\\
$(III)$. Same as in part $(II).$ \qed \\

The following result shows that the characteristic polynomial of a sidigraph in $\Delta^1_n$ is of the form $(1)$.\\

\noindent {\bf Theorem 2.5.} If $S\in \Delta^1_n$, then\\

\begin{equation}
\phi_S(z)=z^n+\sum\limits_{j=1}^{\lfloor{\frac{n}{2}}\rfloor}(-1)^j c_{2j}(S)z^{n-2j},
\end{equation}
where $c_{2j}(S)=|\pounds_{2j}|$ is the cardinality of the set $\pounds_{2j}$.\\
\noindent {\bf Proof.} Let $\phi_S(z)=z^n+\sum\limits_{j=1}^{n}b_j(S)z^{n-j}$. By Theorem $1.1$, we have $b_j(S)=\sum\limits_{L\in\pounds_j}(-1)^{p(L)}s(L)$, where $s(L)=\prod\limits_{Z\in c(L)} s(Z)$. By Lemma $2.4$, for all $j=1,2,\cdots$, we have $b_{2j-1}(S)=0.$\\
Also,
\begin{equation*}
\begin{split}
b_{4j}(S)&=\sum\limits_{L\in\pounds_{4j}}(-1)^{p(L)}s(L)\\&=\sum\limits_{L\in\pounds^1_{4j}}(-1)^{p(L)}s(L)+\sum\limits_{L\in\pounds^2_{4j}}(-1)^{p(L)}s(L)+\sum\limits_{L\in\pounds^3_{4j}}(-1)^{p(L)}s(L),
\end{split}
\end{equation*}

where $\sum\limits_{L\in\pounds^1_{4j}}$ denotes the sum over those linear subsidigraphs $L\in\pounds _{4j}$ whose components are cycles of length $\equiv 0\pmod 4$ only, $\sum\limits_{l\in\pounds^2_{4j}}$
denotes the sum over those linear subsidigraphs $L\in\pounds _{4j}$ whose components are cycles of length $\equiv 2\pmod 4$ only and $\sum\limits_{l\in\pounds^3_{4j}}$ denotes the sum over those linear subsidigraphs $L\in\pounds _{4j}$ which have components consisting of both types of cycles.\\
Now $$\sum\limits_{L\in\pounds^1_{4j}}(-1)^{p(L)}s(L)=\sum\limits_{I}(-1)^{p(L)}s(L)+\sum\limits_{II}(-1)^{p(L)}s(L),$$\\
where $\sum\limits_{I}$ denotes sum over those $L\in\pounds^1_{4j}$ which have an even number of cycles of length  $\equiv 0\pmod 4$ and $\sum\limits_{II}$ denotes sum over those $L\in\pounds^1_{4j}$ which have an odd number of cycles of length  $\equiv 0\pmod 4$. \\
Therefore, by Lemma $2.4$, we have
\begin{equation*}
\begin{split}
\sum\limits_{L\in\pounds^1_{4j}}(-1)^{p(L)}s(L)&=\sum\limits_{I}(-1)^{even}(+1)+\sum\limits_{II}(-1)^{odd}(-1)\\&=\sum\limits_{I}1+\sum\limits_{II}1=|\pounds^1_{4j}|.
\end{split}
\end{equation*}
Now, by Lemma $2.4$, we have\\
\begin{equation*}
\begin{split}
\sum\limits_{L\in\pounds^2_{4j}}(-1)^{p(L)}s(L)&=\sum\limits_{L\in\pounds^2_{4j}}(-1)^{even}(+1)\\&=|\pounds^2_{4j}|.
\end{split}
\end{equation*}

Again by Lemma $2.4$, we have\\
$$\sum\limits_{L\in\pounds^3_{4j}}(-1)^{p(L)}s(L)=\sum\limits_{I}(-1)^{p(L)}s(L)+\sum\limits_{II}(-1)^{p(L)}s(L),$$\\
where $\sum\limits_{I}$ denotes sum over those $L\in\pounds^3_{4j}$ which have an even number of cycles of length  $\equiv 0\pmod 4$ and $\sum\limits_{II}$ denotes sum over those $L\in\pounds^3_{4j}$ which have an odd number of cycles of length  $\equiv 0\pmod 4$. Note that the number of cycles of length $\equiv 2\pmod 4$ is even. \\
Therefore,
\begin{equation*}
\begin{split}
\sum\limits_{L\in\pounds^3_{4j}}(-1)^{p(L)}s(L)&=\sum\limits_{I}(-1)^{even}(+1)+\sum\limits_{II}(-1)^{odd}(-1)\\&=\sum\limits_{I}1+\sum\limits_{II}1=|\pounds^3_{4j}|.
\end{split}
\end{equation*}
Thus $b_{4j}(S)=|\pounds^1_{4j}|+|\pounds^2_{4j}|+|\pounds^3_{4j}|=|\pounds_{4j}|.$\\
Also,
\begin{equation*}
\begin{split}
b_{4j+2}(S)&=\sum\limits_{L\in\pounds_{4j+2}}(-1)^{p(L)}s(L)\\&=\sum\limits_{l\in\pounds^1_{4j+2}}(-1)^{p(L)}s(L)+\sum\limits_{L\in\pounds^2_{4j+2}}(-1)^{p(L)}s(L),
\end{split}
\end{equation*}
where $\sum\limits_{l\in\pounds^1_{4j+2}}$ denotes the sum over those linear subsidigraphs $L\in\pounds _{4j+2}$ whose components are cycles of length $\equiv 2\pmod 4$ only, and $\sum\limits_{l\in\pounds^2_{4j+2}}$ denotes the sum over those linear subsidigraphs $L\in\pounds _{4j+2}$ which have components consisting of both types of cycles.\\
By Lemma $2.4$, we have\\
$$\sum\limits_{L\in\pounds^1_{4j+2}}(-1)^{p(L)}s(L)=\sum\limits_{L\in \pounds^1_{4j+2}}(-1)^{odd}(+1)=-|\pounds^1_{4j+2}|.$$\\
Also,
$$\sum\limits_{L\in\pounds^2_{4j+2}}(-1)^{p(L)}s(L)=\sum\limits_{I}(-1)^{p(L)}s(L)+\sum\limits_{II}(-1)^{p(L)}s(L),$$\\
where $\sum\limits_{I}$ denotes sum over those $L\in\pounds^2_{4j+2}$ which have an even number of cycles of length  $\equiv 0 \pmod 4$ and $\sum\limits_{II}$ denotes sum over those $L\in\pounds^2_{4j+2}$ which have an odd number of cycles of length $\equiv 0\pmod 4$.\\

Again, by Lemma $2.4$, we have\\
\begin{equation*}
\begin{split}
\sum\limits_{L\in\pounds^2_{4j+2}}(-1)^{p(L)}s(L)&=\sum\limits_{I}(-1)^{odd}(+1)+\sum\limits_{II}(-1)^{even}(-1)\\&=\sum\limits_{I}(-1)+\sum\limits_{II}(-1)=-|\pounds^2_{4j+2}|.
\end{split}
\end{equation*}
Therefore, $b_{4j+2}(S)=-|\pounds^1_{4j+2}|-|\pounds^2_{4j+2}|=-|\pounds_{4j+2}|.$\\
Thus we conclude that\\
$$\phi_S(z)=z^n+\sum\limits_{j=1}^{\lfloor{\frac{n}{2}}\rfloor}(-1)^j c_{2j}(S)z^{n-2j},$$\\
where $c_{2j}(S)=|\pounds_{2j}|$ is the cardinality of the set $\pounds_{2j}$. \qed \\

\noindent {\bf Remark 2.6.} Here we note that there exist bipartite and non bipartite non cycle balanced sidigraphs not in $\Delta^1_n$ which have characteristic polynomial with alternating coefficients. Sidigraphs $S_1$ and $S_2$ in Fig. $2$ clearly does not belong to $\Delta^1_n$. By Theorem $1.1$, $\phi_{S_1}(z)=z^6-z^4+2z^2$ and $\phi_{S_2}(z)=z^6-1$.\\

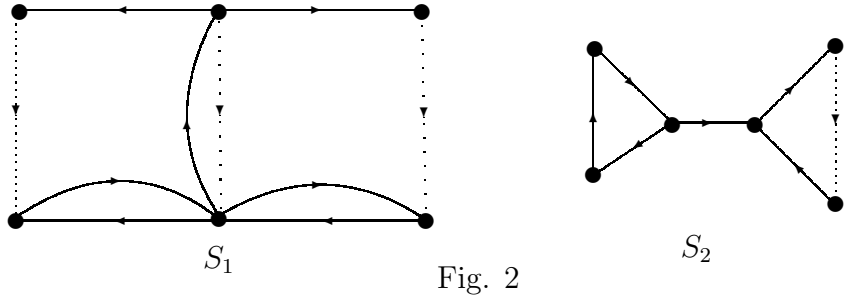
\begin{align*}
%TeXCAD Picture [f1.pic]. Options:
%\grade{\on}
%\emlines{\off}
%\epic{\off}
%\beziermacro{\on}
%\reduce{\on}
%\snapping{\off}
%\quality{8.00}
%\graddiff{0.01}
%\snapasp{1}
%\zoom{4.0000}
\unitlength 1mm % = 2.85pt
\linethickness{0.4pt}
\ifx\plotpoint\undefined\newsavebox{\plotpoint}\fi % GNUPLOT compatibility
\begin{picture}(133.87,40.31)(0,0)
\put(24.25,39.25){\circle*{2}}
\put(50.75,39.25){\circle*{2}}
\put(77.75,39.25){\circle*{2}}
%\vector[middle](50.5,39.5)(23.75,39.5)
\put(37.13,39.5){\vector(-1,0){.07}}\put(50.5,39.5){\line(-1,0){26.75}}
%\end
%\vector[middle](50.75,39.5)(78,39.5)
\put(64.38,39.5){\vector(1,0){.07}}\put(50.75,39.5){\line(1,0){27.25}}
%\end
%\vector[middle](100.5,17.25)(100.5,34.75)
\put(100.5,26){\vector(0,1){.07}}\put(100.5,17.25){\line(0,1){17.5}}
%\end
%\vector[middle](100.75,34.5)(111.25,24.25)
\put(106,29.38){\vector(1,-1){.07}}\multiput(100.75,34.5)(.034539474,-.033717105){304}{\line(1,0){.034539474}}
%\end
%\vector[middle](121.75,24)(132.75,35)
\put(127.25,29.5){\vector(1,1){.07}}\multiput(121.75,24)(.033639144,.033639144){327}{\line(0,1){.033639144}}
%\end
%\vector[middle](132.5,13.75)(122,24.25)
\put(127.25,19){\vector(-1,1){.07}}\multiput(132.5,13.75)(-.033653846,.033653846){312}{\line(0,1){.033653846}}
%\end
%\vector[middle]{dot}(132.75,35)(132.75,13.5)
\put(132.75,24.25){\vector(0,-1){.07}}\multiput(132.68,34.93)(0,-.9773){23}{{\rule{.4pt}{.4pt}}}
%\end
\put(111,24.25){\circle*{2}}
\put(122,24.25){\circle*{2}}
\put(132.75,34.75){\circle*{2}}
\put(132.75,13.75){\circle*{2}}
%\vector[middle]{dot}(23.75,39.5)(23.75,11.5)
\put(23.75,25.5){\vector(0,-1){.07}}\multiput(23.68,39.43)(0,-.9655){30}{{\rule{.4pt}{.4pt}}}
%\end
\put(23.75,11.5){\circle*{2}}
%\qbezvec[middle](51,11.75)(42.25,24.5)(50.5,39.25)
\put(46.5,25){\vector(0,1){.07}}\qbezier(51,11.75)(42.25,24.5)(50.5,39.25)
%\end
%\qbezvec[middle](23.5,11.75)(38.25,21.88)(51,11.5)
\put(37.75,16.75){\vector(1,0){.07}}\qbezier(23.5,11.75)(38.25,21.88)(51,11.5)
%\end
%\vector[middle]{dot}{2}(50.75,39.75)(51,11.25)
\put(50.88,25.5){\vector(0,-1){.07}}\multiput(50.68,39.68)(.0167,-1.9){16}{{\rule{.4pt}{.4pt}}}
%\end
%\vector[middle](50.25,11.5)(23.5,11.5)
\put(36.88,11.5){\vector(-1,0){.07}}\put(50.25,11.5){\line(-1,0){26.75}}
%\end
%\vector[middle](78,11.5)(51,11.5)
\put(64.5,11.5){\vector(-1,0){.07}}\put(78,11.5){\line(-1,0){27}}
%\end
%\vector[middle]{dot}{2}(77.75,39.25)(78.25,11.25)
\put(78,25.25){\vector(0,-1){.07}}\multiput(77.68,39.18)(.0333,-1.8667){16}{{\rule{.4pt}{.4pt}}}
%\end
\put(50.75,11.75){\circle*{2}}
\put(78.25,11.5){\circle*{2}}
%\qbezvec[middle](50.5,12)(64.75,20.63)(78,11.75)
\put(64.5,16.25){\vector(1,0){.07}}\qbezier(50.5,12)(64.75,20.63)(78,11.75)
%\end
\put(50.75,6.25){\makebox(0,0)[cc]{$S_1$}}
\put(114.25,7.75){\makebox(0,0)[cc]{$S_2$}}
\put(85.25,3.5){\makebox(0,0)[cc]{Fig. $2$}}
\put(100.75,34.25){\circle*{2.06}}
\put(100.5,17.5){\circle*{2}}
%\vector[middle](111,24.5)(100.5,17.5)
\put(105.75,21){\vector(-3,-2){.07}}\multiput(111,24.5)(-.05048077,-.03365385){208}{\line(-1,0){.05048077}}
%\end
%\vector[middle](110.75,24.5)(122,24.5)
\put(116.38,24.5){\vector(1,0){.07}}\put(110.75,24.5){\line(1,0){11.25}}
%\end
\end{picture}
\end{align*}

The following result shows that the characteristic polynomial of a sidigraph in $\Delta^2_n$ is of the form $(2)$. Proof is same as the proof of Theorem $2.5$.\\

\noindent {\bf Theorem 2.7.} Let $S\in \Delta^2_n$. Then characteristic polynomial is given by\\
$$\phi_S(z)=z^n+\sum\limits_{j=1}^{\lfloor{\frac{n}{2}}\rfloor} c_{2j}(S)z^{n-2j},$$\\
where $c_{2j}(S)=|\pounds_{2j}|$ is the cardinality of the set $\pounds_{2j}$.\\

%TeXCAD Picture [f2.pic]. Options:
%\grade{\on}
%\emlines{\off}
%\epic{\off}
%\beziermacro{\on}
%\reduce{\on}
%\snapping{\off}
%\quality{8.00}
%\graddiff{0.01}
%\snapasp{1}
%\zoom{4.0000}
\unitlength 1mm % = 2.85pt
\linethickness{0.4pt}
\ifx\plotpoint\undefined\newsavebox{\plotpoint}\fi % GNUPLOT compatibility
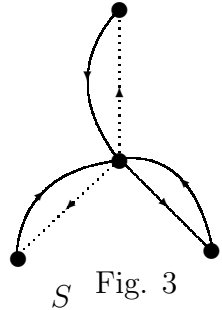
\begin{picture}(51.56,46.29)(0,0)
\put(38.25,25.5){\circle*{2}}
%\vector[middle](38.5,25.5)(51,13)
\put(44.75,19.25){\vector(1,-1){.07}}\multiput(38.5,25.5)(.033692722,-.033692722){371}{\line(0,-1){.033692722}}
%\end
%\vector[middle]{dot}(38.25,25.75)(38.25,45)
\put(38.25,35.38){\vector(0,1){.07}}\multiput(38.18,25.68)(0,.9625){21}{{\rule{.4pt}{.4pt}}}
%\end
\put(38.25,45.5){\circle*{2}}
\put(50.5,13.5){\circle*{2}}
\put(24.75,12.5){\circle*{2}}
%\qbezvec[middle](38.25,45.75)(29.75,36.38)(38.25,25.5)
\put(34,36){\vector(0,-1){.07}}\qbezier(38.25,45.75)(29.75,36.38)(38.25,25.5)
%\end
%\qbezvec[middle](24.5,12.5)(25.25,24)(38,25.5)
\put(28.25,21.5){\vector(1,1){.07}}\qbezier(24.5,12.5)(25.25,24)(38,25.5)
%\end
%\qbezvec[middle](50.75,13.25)(48,27)(38.25,25.75)
\put(46.25,23.25){\vector(-1,1){.07}}\qbezier(50.75,13.25)(48,27)(38.25,25.75)
%\end
%\vector[middle]{dot}(38,25.5)(24.25,12.5)
\put(31.13,19){\vector(-1,-1){.07}}\multiput(37.93,25.43)(-.6548,-.619){22}{{\rule{.4pt}{.4pt}}}
%\end
\put(30.5,7.5){\makebox(0,0)[cc]{$S$}}
\put(40.5,9){\makebox(0,0)[cc]{Fig. $3$}}
\end{picture}

\noindent {\bf Remark 2.8.} We note that there exist bipartite and non bipartite non cycle balanced sidigraphs not in $\Delta^2_n$ which have characteristic polynomial of the form $(2)$. Sidigraphs $S_1$ and $S_2$ in Fig. $1$ and sidigraph $S$ in Fig. $3$ does not belong to  $\Delta^2_n$  because former are non bipartite and latter has a positive cycle of length $2$. By Theorem $1.1$, $\phi_{S_1}(z)=z^{17}+3z^{11}+z^5$, $\phi_{S_2}(z)=z^{17}+z^{11}+z^5$ and $\phi_S(z)=z^4+z^2.$\\

\indent Recall the definition of Cartesian product of two sidigraphs from \cite{pm}. Let $S_1=(V_1,\mathscr{A}_1,\sigma_1)$ and $S_2=(V_2,\mathscr{A}_2,\sigma_2)$ be two sidigraphs, their Cartesian product (or sum) denoted by $S_1\times S_2$ is the sidigraph $(V_1\times V_2,\mathscr{A},\sigma)$, where the arc set $\mathscr{A}$ is that of the Cartesian product of underlying unsigned digraphs and the sign function is defined by\\
 $$\sigma((u_i,v_j),(u_k,v_l))=\left \{\begin{array}{lr}\sigma_1(u_i,u_k), &\mbox{if $j=l$},\\
 \sigma_2(v_j,v_l), &\mbox{if $i=k.$}
 \end{array} \right.$$\\

Unlike Kronecker product \cite{ma}, Cartesian product of two  strongly connected sidigraphs is always strongly connected as can be seen in the following result.\\

\noindent {\bf Lemma 2.9.} Let $S_1$ and $S_2$ be two strongly connected sidigraphs. Then $S_1\times S_2$ is strongly connected.\\
\noindent {\bf Proof.} Let $(u_i,v_j),(u_p,v_q)\in V(S_1\times S_2)$, where we assume $p\le q$ (case $p > q$ can be dealt similarly). Since $S_1$ is strongly connected, there exists a directed path $(u_i,u_{i+1})(u_{i+1},u_{i+2})\cdots(u_{p-1},u_p)$. Also, strong connectedness of $S_2$ implies there exists a directed path  $(v_j,v_{j+1})(v_{j+1},v_{j+2})\cdots(v_{q-1},v_q)$. By definition of Cartesian product, Fig. $4$ illustrates that there exists a directed path from $(u_i,v_j)$ to $(u_p,v_q)$. Signs do not play any role in connectedness, so we take all arcs in Fig. $4$ positive. Similarly, one can prove the reverse part. \qed \\

%TeXCAD Picture [s1.pic]. Options:
%\grade{\on}
%\emlines{\off}
%\epic{\off}
%\beziermacro{\on}
%\reduce{\on}
%\snapping{\off}
%\quality{8.00}
%\graddiff{0.01}
%\snapasp{1}
%\zoom{4.0000}
\unitlength 1mm % = 2.85pt
\linethickness{0.4pt}
\ifx\plotpoint\undefined\newsavebox{\plotpoint}\fi % GNUPLOT compatibility
\begin{picture}(124.25,36.25)(0,0)
\put(8.75,30.25){\circle*{2}}
\put(16,20){\circle*{2}}
\put(26.5,30.25){\circle*{2}}
%\vector[middle](8.5,30.5)(15.75,19.75)
\put(12.13,25.13){\vector(2,-3){.07}}\multiput(8.5,30.5)(.03372093,-.05){215}{\line(0,-1){.05}}
%\end
%\vector[middle](16,20)(26,30.25)
\put(21,25.13){\vector(1,1){.07}}\multiput(16,20)(.033670034,.034511785){297}{\line(0,1){.034511785}}
%\end
%\vector[middle](26.5,30.25)(38,18.5)
\put(32.25,24.38){\vector(1,-1){.07}}\multiput(26.5,30.25)(.03372434,-.034457478){341}{\line(0,-1){.034457478}}
%\end
\put(37.5,18.75){\circle*{2}}
~~\put(41.5,7.25){\makebox(0,0)[cc]{Fig. $4$}}
\put(15.75,15){\makebox(0,0)[cc]{$(u_i,v_{j+1})$}}
\put(27.25,34.75){\makebox(0,0)[cc]{$(u_{i+1},v_{j+1})$}}
\put(37.25,14.5){\makebox(0,0)[cc]{$(u_{i+1},v_{j+2})$}}
\put(7.75,35){\makebox(0,0)[cc]{$(u_i,v_j)$}}
\put(51,21.5){\makebox(0,0)[cc]{$\cdots$}}
\put(63.25,18){\circle*{2}}
\put(69.25,31){\circle*{2}}
\put(84.75,30.75){\circle*{2}}
%\vector[middle](63,18.25)(69,31)
\put(66,24.63){\vector(1,2){.07}}\multiput(63,18.25)(.03370787,.07162921){178}{\line(0,1){.07162921}}
%\end
%\vector[middle](69,31)(84.75,31)
\put(76.88,31){\vector(1,0){.07}}\put(69,31){\line(1,0){15.75}}
%\end
\put(62,14){\makebox(0,0)[cc]{$(u_{p-1},v_p)$}}
\put(68.25,36){\makebox(0,0)[cc]{$(u_p,v_p)$}}
\put(85.25,36.25){\makebox(0,0)[cc]{$(u_p,v_{p+1})$}}
\put(100.25,22.25){\makebox(0,0)[cc]{$\cdots$}}
\put(112.25,30.75){\circle*{2}}
\put(112.75,36.25){\makebox(0,0)[cc]{$(u_p,v_{q-1})$}}
\put(123.25,30.5){\circle*{2}}
%\vector[middle](114,31)(122,31)
\put(118,31){\vector(1,0){.07}}\put(114,31){\line(1,0){8}}
%\end
\put(123,25.25){\makebox(0,0)[cc]{$(u_p,v_q)$}}
\end{picture}

\noindent {\bf Definition 2.10.} A sidigraph is said to integral (real or Gaussian) according as spectrum of $S$ is integral (real or Gaussian) respectively.\\

Integral, real and Gaussian spectral digraphs were studied by Esser and Harary \cite{eh1}. The following three results show the existence of non cycle balanced integral, real and Gaussian sidigraphs.\\

%TeXCAD Picture [f3.pic]. Options:
%\grade{\on}
%\emlines{\off}
%\epic{\off}
%\beziermacro{\on}
%\reduce{\on}
%\snapping{\off}
%\quality{8.00}
%\graddiff{0.01}
%\snapasp{1}
%\zoom{4.0000}
\unitlength 1mm % = 2.85pt
\linethickness{0.4pt}
\ifx\plotpoint\undefined\newsavebox{\plotpoint}\fi % GNUPLOT compatibility
\begin{picture}(86.75,58.75)(0,0)
\put(16.75,48.5){\circle*{2}}
\put(44.75,48.75){\circle*{2}}
\put(58,48.5){\circle*{2}}
\put(85.5,48.5){\circle*{2}}
%\vector[middle](17,49)(44.75,49)
\put(30.88,49){\vector(1,0){.07}}\put(17,49){\line(1,0){27.75}}
%\end
%\vector[middle](44.75,48.75)(44.75,22)
\put(44.75,35.38){\vector(0,-1){.07}}\put(44.75,48.75){\line(0,-1){26.75}}
%\end
%\vector[middle](44.75,22.5)(16,22.5)
\put(30.38,22.5){\vector(-1,0){.07}}\put(44.75,22.5){\line(-1,0){28.75}}
%\end
%\vector[middle](16.75,22.75)(16.75,49)
\put(16.75,35.88){\vector(0,1){.07}}\put(16.75,22.75){\line(0,1){26.25}}
%\end
%\vector[middle](58,49)(85.25,49)
\put(71.63,49){\vector(1,0){.07}}\put(58,49){\line(1,0){27.25}}
%\end
%\vector[middle](85,49)(85,22)
\put(85,35.5){\vector(0,-1){.07}}\put(85,49){\line(0,-1){27}}
%\end
%\vector[middle](84.75,22.75)(57.25,22.75)
\put(71,22.75){\vector(-1,0){.07}}\put(84.75,22.75){\line(-1,0){27.5}}
%\end
%\vector[middle](58,22.75)(58,49)
\put(58,35.88){\vector(0,1){.07}}\put(58,22.75){\line(0,1){26.25}}
%\end
\put(16.5,22.75){\circle*{2}}
\put(44.75,22.5){\circle*{2}}
\put(58,22.75){\circle*{2}}
\put(85,22.5){\circle*{2}}
%\vector[middle]{dot}(16.5,23)(45,49)
\put(30.75,36){\vector(1,1){.07}}\multiput(16.43,22.93)(.69512,.63415){42}{{\rule{.4pt}{.4pt}}}
%\end
%\vector[middle]{dot}(85.25,49)(58,23.25)
\put(71.63,36.13){\vector(-1,-1){.07}}\multiput(85.18,48.93)(-.69872,-.66026){40}{{\rule{.4pt}{.4pt}}}
%\end
%\vector[middle]{dot}(58.25,49)(84.75,23)
\put(71.5,36){\vector(1,-1){.07}}\multiput(58.18,48.93)(.6625,-.65){41}{{\rule{.4pt}{.4pt}}}
%\end
%\qbezvec[middle](45,49.25)(29.63,58.5)(16.75,48.75)
\put(30.25,53.75){\vector(-1,0){.07}}\qbezier(45,49.25)(29.63,58.5)(16.75,48.75)
%\end
%\qbezvec[middle](16.75,49)(8.38,36.5)(16.5,23)
\put(12.5,36.25){\vector(0,-1){.07}}\qbezier(16.75,49)(8.38,36.5)(16.5,23)
%\end
%\qbezvec[middle](16.5,23)(30.38,12.5)(44.75,23)
\put(30.5,17.75){\vector(1,0){.07}}\qbezier(16.5,23)(30.38,12.5)(44.75,23)
%\end
%\qbezvec[middle](85.25,48.75)(70.88,58.75)(58,48.75)
\put(71.25,53.75){\vector(-1,0){.07}}\qbezier(85.25,48.75)(70.88,58.75)(58,48.75)
%\end
%\qbezvec[middle](58,49.25)(48.5,35.75)(58,23.25)
\put(53.25,36){\vector(0,-1){.07}}\qbezier(58,49.25)(48.5,35.75)(58,23.25)
%\end
%\qbezvec[middle](57.75,23)(69.75,14)(84.75,23)
\put(70.5,18.5){\vector(1,0){.07}}\qbezier(57.75,23)(69.75,14)(84.75,23)
%\end
\put(24.5,12.25){\makebox(0,0)[cc]{$S_1$}}
\put(70,12.25){\makebox(0,0)[cc]{$S_2$}}
\put(50,11.75){\makebox(0,0)[cc]{Fig. $5$}}
\end{picture}

\noindent {\bf Theorem 2.11.} For each positive integer $n\ge 4$, there exists a family of $n$ integral cospectral, strongly connected, non symmetric and non cycle balanced sidigraphs of order $4^n$.\\
\noindent{\bf Proof.} Consider sidigraphs $S_1$ and $S_2$ in Fig. $5$. Clearly $S_1$ and $S_2$ are non cycle balanced and strongly connected. By Theorem $1.1$, \\
$$\phi_{S_1}(z)=\phi_{S_2}(z)=z^4-3z^2+2z.$$\\
Therefore, $spec(S_1)=spec(S_2)= \{-2,0,1^{(2)}\}$. That is, $S_1$ and $S_2$ are integral cospectral.\\
Let $$S^{(k)}=S_1\times S_1\times\cdots\times S_1\times S_2\times S_2\times\cdots\times S_2,$$\\
where we take $k$ copies of $S_1$ and $n-k$ copies of $S_2$. Clearly, for each $n$, we have $n$ cospectral sidigraphs $S^{(k)},~~k=1,2,\cdots,n$ of order $4^n$. $S_1$ and $S_2$ are non symmetric implies $S^{(k)}$ is non symmetric. By repeated application of Lemma $2.9$ and using the fact that the Cartesian product of sidigraphs is cycle balanced if and only if the constituent sidigraphs are cycle balanced [Theorem $4.8$, \cite {pm}], the result follows.\qed \\

Integral sidigraphs are obviously real. There exists non integral real sidigraphs as can be see in the following result.\\

%TeXCAD Picture [f4.pic]. Options:
%\grade{\on}
%\emlines{\off}
%\epic{\off}
%\beziermacro{\on}
%\reduce{\on}
%\snapping{\off}
%\quality{8.00}
%\graddiff{0.01}
%\snapasp{1}
%\zoom{4.0000}
\unitlength 1mm % = 2.85pt
\linethickness{0.4pt}
\ifx\plotpoint\undefined\newsavebox{\plotpoint}\fi % GNUPLOT compatibility
\begin{picture}(116.56,63.38)(0,0)
\put(13.25,52.25){\circle*{2}}
\put(40.75,52.25){\circle*{2}}
\put(13.25,19.75){\circle*{2}}
\put(40.75,19.75){\circle*{2}}
\put(52,52){\circle*{2}}
\put(52,19.5){\circle*{2}}
\put(79.5,52){\circle*{2}}
\put(87.5,51.75){\circle*{2}}
\put(115.5,51.75){\circle*{2}}
\put(115.5,19.5){\circle*{2}}
\put(87.5,19.25){\circle*{2}}
%\vector[middle]{dot}(40.5,52.5)(40.75,20.25)
\put(40.63,36.38){\vector(0,-1){.07}}\multiput(40.43,52.43)(.00758,-.97727){34}{{\rule{.4pt}{.4pt}}}
%\end
%\vector[middle]{dot}(40.75,20.25)(40.75,20.25)
\put(40.75,20.25){\vector(0,1){.07}}\multiput(40.68,20.18)(0,0){3}{{\rule{.4pt}{.4pt}}}
%\end
%\vector[middle]{dot}(79.25,52.25)(79.5,21.5)
\put(79.38,36.88){\vector(0,-1){.07}}\multiput(79.18,52.18)(.0081,-.9919){32}{{\rule{.4pt}{.4pt}}}
%\end
%\vector[middle]{dot}(115.5,52)(115.5,18.5)
\put(115.5,35.25){\vector(0,-1){.07}}\multiput(115.43,51.93)(0,-.98529){35}{{\rule{.4pt}{.4pt}}}
%\end
%\vector[middle](13,20.25)(13,52.5)
\put(13,36.38){\vector(0,1){.07}}\put(13,20.25){\line(0,1){32.25}}
%\end
%\vector[middle](13.25,52.5)(40.75,52.5)
\put(27,52.5){\vector(1,0){.07}}\put(13.25,52.5){\line(1,0){27.5}}
%\end
%\vector[middle](40.75,20)(13,20)
\put(26.88,20){\vector(-1,0){.07}}\put(40.75,20){\line(-1,0){27.75}}
%\end
%\vector[middle](51.75,52.25)(79.75,52.25)
\put(65.75,52.25){\vector(1,0){.07}}\put(51.75,52.25){\line(1,0){28}}
%\end
%\vector[middle](51.75,19.75)(51.75,52.5)
\put(51.75,36.13){\vector(0,1){.07}}\put(51.75,19.75){\line(0,1){32.75}}
%\end
%\vector[middle](79.25,20)(51.25,20)
\put(65.25,20){\vector(-1,0){.07}}\put(79.25,20){\line(-1,0){28}}
%\end
%\vector[middle](87.5,19.25)(87.5,52.5)
\put(87.5,35.88){\vector(0,1){.07}}\put(87.5,19.25){\line(0,1){33.25}}
%\end
%\vector[middle](87.75,52.5)(116,52.5)
\put(101.88,52.5){\vector(1,0){.07}}\put(87.75,52.5){\line(1,0){28.25}}
%\end
%\vector[middle](115.5,19.75)(87,19.75)
\put(101.25,19.75){\vector(-1,0){.07}}\put(115.5,19.75){\line(-1,0){28.5}}
%\end
%\vector[middle](13.25,20.5)(41,52.5)
\put(27.13,36.5){\vector(3,4){.07}}\multiput(13.25,20.5)(.0337181045,.0388821385){823}{\line(0,1){.0388821385}}
%\end
%\qbezvec[middle](40.75,52.75)(27.63,62.75)(13,52.75)
\put(27.25,57.75){\vector(-1,0){.07}}\qbezier(40.75,52.75)(27.63,62.75)(13,52.75)
%\end
%\qbezvec[middle](13,52.5)(3,37.75)(13,20)
\put(8,37){\vector(0,-1){.07}}\qbezier(13,52.5)(3,37.75)(13,20)
%\end
%\qbezvec[middle](13,20.5)(23.25,9.75)(40.5,20)
\put(25,15){\vector(1,0){.07}}\qbezier(13,20.5)(23.25,9.75)(40.5,20)
%\end
%\qbezvec[middle](79.5,52.25)(67.38,62.75)(51.75,52.25)
\put(66.5,57.5){\vector(-1,0){.07}}\qbezier(79.5,52.25)(67.38,62.75)(51.75,52.25)
%\end
%\qbezvec[middle](51.75,52)(41.75,36)(51.75,20)
\put(46.75,36){\vector(0,-1){.07}}\qbezier(51.75,52)(41.75,36)(51.75,20)
%\end
%\qbezvec[middle](51.5,20.25)(64.25,9.75)(79,20.25)
\put(64.75,15){\vector(1,0){.07}}\qbezier(51.5,20.25)(64.25,9.75)(79,20.25)
%\end
%\qbezvec[middle](115.75,52.25)(100.75,63.38)(87.75,52)
\put(101.25,57.75){\vector(-1,0){.07}}\qbezier(115.75,52.25)(100.75,63.38)(87.75,52)
%\end
%\qbezvec[middle](87.75,52.25)(78.88,38.38)(87.5,20)
\put(83.25,37.25){\vector(0,-1){.07}}\qbezier(87.75,52.25)(78.88,38.38)(87.5,20)
%\end
%\qbezvec[middle](87.25,19.75)(99.63,9.13)(115.5,20)
\put(100.5,14.5){\vector(1,0){.07}}\qbezier(87.25,19.75)(99.63,9.13)(115.5,20)
%\end
\put(20.75,8.25){\makebox(0,0)[cc]{$S_1$}}
\put(99,9.5){\makebox(0,0)[cc]{$S_3$}}
\put(63,9.25){\makebox(0,0)[cc]{$S_2$}}
\put(57.25,2.5){\makebox(0,0)[cc]{Fig. $6$}}
%\vector[middle](79.75,52.25)(51.75,20)
\put(65.75,36.13){\vector(-1,-1){.07}}\multiput(79.75,52.25)(-.0337349398,-.0388554217){830}{\line(0,-1){.0388554217}}
%\end
%\vector[middle]{dot}(52,52.25)(79,20.25)
\put(65.5,36.25){\vector(3,-4){.07}}\multiput(51.93,52.18)(.61364,-.72727){45}{{\rule{.4pt}{.4pt}}}
%\end
\put(78.75,20.25){\circle*{2}}
\end{picture}

\noindent {\bf Theorem 2.12.}  For each positive integer $n\ge 4$, there exists a family of $n$ real cospectral, strongly connected, non symmetric and non cycle balanced sidigraphs of order $4^n$.\\
\noindent{\bf Proof.} Consider sidigraphs $S_1$, $S_2$ and $S_3$ shown in Fig. $6$. Clearly, all three sidigraphs are non cycle balanced and strongly connected. By Theorem $1.1$,\\
$$\phi_{S_1}(z)=\phi_{S_2}(z)=\phi_{S_3}(z)=z^4-3z^2+2.$$\\
Therefore, $spec(S_1)=spec(S_2)=spec(S_3)=\{-\sqrt{2},-1,1,\sqrt{2}\}$. Take any two sidigraphs among $S_1$, $S_2$ and $S_3$ and apply procedure of Theorem $2.11$, the result follows.\qed \\

\begin{align*}
%TeXCAD Picture [s3.pic]. Options:
%\grade{\on}
%\emlines{\off}
%\epic{\off}
%\beziermacro{\on}
%\reduce{\on}
%\snapping{\off}
%\quality{8.00}
%\graddiff{0.01}
%\snapasp{1}
%\zoom{4.0000}
\unitlength 1mm % = 2.85pt
\linethickness{0.4pt}
\ifx\plotpoint\undefined\newsavebox{\plotpoint}\fi % GNUPLOT compatibility
\begin{picture}(90.06,40.78)(0,0)
\put(8.5,39.25){\circle*{2}}
\put(28.75,39.75){\circle*{2}}
\put(28.25,22){\circle*{2}}
\put(8.5,22.25){\circle*{2}}
\put(66.25,37.75){\circle*{2}}
\put(88.5,38){\circle*{2}}
\put(89,22){\circle*{2}}
\put(66,22.25){\circle*{2}}
%\vector[middle]{dot}(8.25,22.5)(28.25,22)
\put(18.25,22.25){\vector(1,0){.07}}\multiput(8.18,22.43)(.9524,-.0238){22}{{\rule{.4pt}{.4pt}}}
%\end
%\vector[middle](28.75,39.75)(8,39.75)
\put(18.38,39.75){\vector(-1,0){.07}}\put(28.75,39.75){\line(-1,0){20.75}}
%\end
%\vector[middle](8,39.75)(8,22.5)
\put(8,31.13){\vector(0,-1){.07}}\put(8,39.75){\line(0,-1){17.25}}
%\end
%\vector[middle]{dot}(88.5,38)(65.75,38.25)
\put(77.13,38.13){\vector(-1,0){.07}}\multiput(88.43,37.93)(-.9891,.0109){24}{{\rule{.4pt}{.4pt}}}
%\end
%\vector[middle]{dot}(28.25,22.5)(28.25,40)
\put(28.25,31.25){\vector(0,1){.07}}\multiput(28.18,22.43)(0,.9722){19}{{\rule{.4pt}{.4pt}}}
%\end
%\vector[middle]{dot}(89,22)(89,38.5)
\put(89,30.25){\vector(0,1){.07}}\multiput(88.93,21.93)(0,.9706){18}{{\rule{.4pt}{.4pt}}}
%\end
%\vector[middle](60.5,21.5)(60.5,38.5)
\put(60.5,30){\vector(0,1){.07}}\put(60.5,21.5){\line(0,1){17}}
%\end
%\vector[middle](66,38)(66,22.5)
\put(66,30.25){\vector(0,-1){.07}}\put(66,38){\line(0,-1){15.5}}
%\end
%\vector[middle](66,22.5)(89,22.5)
\put(77.5,22.5){\vector(1,0){.07}}\put(66,22.5){\line(1,0){23}}
%\end
%\vector[middle](66,23)(88.75,38.25)
\put(77.38,30.63){\vector(3,2){.07}}\multiput(66,23)(.050331858,.033738938){452}{\line(1,0){.050331858}}
%\end
%\vector[middle](89,22.5)(66,38.25)
\put(77.5,30.38){\vector(-3,2){.07}}\multiput(89,22.5)(-.049250535,.03372591){467}{\line(-1,0){.049250535}}
%\end
\put(60.5,21.25){\circle*{1.5}}
\put(48,16.5){\makebox(0,0)[cc]{$S_2$}}
\put(76,17.75){\makebox(0,0)[cc]{$S_3$}}
\put(16,17.75){\makebox(0,0)[cc]{$S_1$}}
\put(42.75,9){\makebox(0,0)[cc]{Fig. $7$}}
%\vector[middle](60.5,38)(33.75,38)
\put(47.13,38){\vector(-1,0){.07}}\put(60.5,38){\line(-1,0){26.75}}
%\end
%\vector[middle](34.25,38)(34.25,21.5)
\put(34.25,29.75){\vector(0,-1){.07}}\put(34.25,38){\line(0,-1){16.5}}
%\end
%\vector[middle](34.25,21.5)(60.5,21.5)
\put(47.38,21.5){\vector(1,0){.07}}\put(34.25,21.5){\line(1,0){26.25}}
%\end
%\vector[middle](60.25,21.5)(34.5,38)
\put(47.38,29.75){\vector(-3,2){.07}}\multiput(60.25,21.5)(-.05255102,.033673469){490}{\line(-1,0){.05255102}}
%\end
%\vector[middle]{dot}(60.75,37.75)(34.25,22)
\put(47.5,29.88){\vector(-3,-2){.07}}\multiput(60.68,37.68)(-.82813,-.49219){33}{{\rule{.4pt}{.4pt}}}
%\end
\put(34.5,38){\circle*{2}}
\put(60.75,37.75){\circle*{2}}
\put(34.25,21.5){\circle*{2}}
\end{picture}
\end{align*}

Every integral sidigraph is obviously Gaussian. The next result shows that there exists non integral Gaussian sidigraphs i.e., sidigraphs with eigenvalues of the form $a+\iota b$, where $a$ and $b$  are integers with $b\neq 0$ for some $b$.\\

\noindent {\bf Theorem 2.13.} For each positive integer $n\ge 4$, there exists a collection of $n$ Gaussian cospectral, strongly connected, non symmetric and non cycle balanced sidigraphs of order $4^n$.\\
\noindent{\bf Proof.} Consider sidigraphs $S_1$, $S_2$ and $S_3$ in Fig. $7$. It is clear that $S_1$ is cycle balanced whereas $S_2$ and $S_3$ are non cycle balanced. Moreover all three sidigraphs are strongly connected. By Theorem $1.1$,\\
$$\phi_{S_1}(z)=\phi_{S_2}(z)=\phi_{S_3}(z)=z^4-1.$$\\
Therefore, $spec(S_1)=spec(S_2)=spec(S_3)=\{-1,1,-\iota,\iota\}$. Hence $S_1$, $S_2$ and $S_3$ are Gaussian cospectral.
Take any two sidigraphs among $S_1$, $S_2$ and $S_3$ and proceed in a similar way as in Theorem $2.11$, the result follows.\qed \\

Two digraphs $D_1$ and $D_2$ are said to be quasi-cospectral if there exist sidigraphs $S_1$ and $S_2$ on $D_1$ and $D_2$ respectively such that $\phi_{S_1}(z)=\phi_{S_2}(z)$ i.e., $S_1$ and $S_2$ are cospectral. Two cospectral digraphs are quasi-cospectral by Theorem $1.2$, as we can take any two cycle balanced sidigraphs one on each digraph. Two digraphs are said to be strictly quasi-cospectral if they are quasi-cospectral but not cospectral. Two digraphs $D_1$ and $D_2$ are said to be strongly quasi-cospectral if both $D_1$ and $D_2$ are cospectral and there exists non cycle balanced sidigraphs respectively $S_1$ and $S_2$ on them such that $\phi_{S_1}(z)=\phi_{S_2}(z)$. It is clear that if $D_1$ and $D_2$ are strongly quasi-cospectral digraphs, then both should have at least on cycle.  For quasi-cospectral and strongly quasi-cospectral graphs and digraphs see \cite{agp,mk}.\\

\noindent {\bf Definition 2.14.} We say two digraphs $D_1$ and $D_2$ are integral, real and Gaussian strongly quasi-cospectral if both $D_1$ and $D_2$ are respectively integral, real and Gaussian cospectral and there exists non cycle balanced sidigraphs $S_1$ and $S_2$ on them which are respectively integral, real and Gaussian cospectral.\\

The following two results show the existence of integral and real strongly quasi-cospectral digraphs.\\

\noindent {\bf Theorem 2.15.} For each positive integer $n\ge 4$, there exists a family of $n$ integral, strongly connected, non symmetric and strongly quasi-cospectral digraphs of order $4^n$.\\
\noindent{\bf Proof.} Let $D_1$ and $D_2$  respectively be the underlying digraphs of integral sidigraphs $S_1$ and $S_2$ shown in Fig. $5$. Then $D_1$ and $D_2$ are all-positive sidigraphs. By Theorem $1.1$, we have \\
$$\phi_{D_1}(z)=\phi_{D_2}(z)=z^4-3z^2-2z.$$\\
Therefore, $spec(D_1)=spec(D_2)=\{-1^{(2)},0,2\}.$\\
Put $D^{(k)}=D_1\times D_1\times\cdots\times D_1\times D_2\times D_2\times\cdots\times D_2,$ where we take $k$ copies of $D_1$ and $n-k$ copies of $D_2$. In this way, for each $n\ge 4$ we get $n$ cospectral non symmetric and strongly connected  integral digraphs. Thus for any two of these integral cospectral digraphs $D^{(k_1)}$ and  $D^{(k_2)}$ there exists corresponding non cycle balanced sidigraphs $S^{(k_1)}$ and $S^{(k_2)}$ on them which are integral cospectral.\qed \\

The following result shows the existence of real strongly quasi-cospectral digraphs.\\

\noindent {\bf Theorem 2.16.} For each positive integer $n\ge 4$, there exists a collection of $n$ real, strongly connected, non symmetric and  strongly quasi-cospectral digraphs of order $4^n$.\\
\noindent{\bf Proof.} Let $D_1$ and $D_2$ be the underlying digraphs of sidigraphs $S_1$ and $S_2$ as shown in Fig. $6$. It is easy to see that $\phi_{D_1}(z)=\phi_{D_2}(z)=z^4-3z^2-2z$ and $spec(D_1)=spec(D_2)=\{-1^{(2)},0,2\}.$ Also $spec(S_1)=spec(S_2)=\{-\sqrt{2},-1,1,\sqrt{2}\}.$\\
\indent Thus $D_1$ and $D_2$ are real strongly quasi-cospectral. Applying the same technique as in Theorem $2.15$, the result follows.\qed \\

\section{Energy of sidigraphs}

In \cite{rgc}, the authors compared the energies of digraphs in $\Delta_n $ consisting of bipartite digraphs with each cycle of length $\equiv 2\pmod 4$. We now derive integral expressions for sidigraphs in $\Delta^1_n$ and $\Delta^2_n$ and compare energies of sidigraphs in $\Delta^1_n$ by means of quasi-order relation.\\

Given sidigraphs $S_1$ and $S_2$ in $\Delta^1_n$, by Theorem $2.5$, for $i=1,2$, we have\\
$$\phi_{S_i(z)}=z^n+\sum\limits_{j=1}^{\lfloor{\frac{n}{2}}\rfloor}(-1)^j c_{2j}(S_i)z^{n-2j},$$
where $c_{2j}(S_i)$ are non negative integers for all $j=1,2,\cdots,\lfloor{\frac{n}{2}}\rfloor$. If $c_{2j}(S_1)\le c_{2j}(S_2)$ for all $j=1,2,\cdots,\lfloor{\frac{n}{2}}\rfloor$, then we define $S_1\preceq S_2$. If in addition $c_{2j}(S_1) < c_{2j}(S_2)$ for some $j=1,2,\cdots,\lfloor{\frac{n}{2}}\rfloor$, then we write $S_1 \prec S_2$. The following result whose proof is same as the proof of the [Theorem $2.4$, \cite{rgc}] shows energy increases with respect to this quasi-order relation.\\

\noindent {\bf Theorem 3.1.} If $S\in \Delta^1_n$, then\\
$$E(S)=\frac{1}{\pi} \int\limits_{-\infty}^{\infty}\frac{1}{z^2} \log[1+\sum\limits_{j=1}^{\lfloor{\frac{n}{2}}\rfloor} c_{2j}(S) z^{2j}]dz.$$\\
In particular, if $S_1,S_2\in \Delta^1_n$ and $S_1\prec S_2$ then $E(S_1)<E(S_2)$.\\

Next we derive an integral expression for the energy of a sidigraph in $\Delta^2_n$.\\

\noindent {\bf Theorem 3.2.} If $S\in \Delta^2_n$, then\\
$$E(S)=\frac{1}{\pi} \int\limits_{-\infty}^{\infty}\frac{1}{z^2} \log|1+\sum\limits_{j=1}^{\lfloor{\frac{n}{2}}\rfloor} (-1)^j c_{2j}(S) z^{2j}|dz.$$\\
\noindent {\bf Proof.} We know from \cite{pm} that the energy of a sidigraph $S$ satisfies the integral expression\\
$$E(S)=\frac{1}{\pi} \int\limits_{-\infty}^{\infty}\frac{1}{z^2} \log|z^n \phi_S(\frac{\iota}{z})|dz.$$\\
Assume $S\in \Delta^2_n$, then $\phi_{S(z)}=z^n+\sum\limits_{j=1}^{\lfloor{\frac{n}{2}}\rfloor} c_{2j}(S)z^{n-2j}$ \\
so that
\begin{align*}
 E(S)&=\frac{1}{\pi} \int\limits_{-\infty}^{\infty}\frac{1}{z^2} \log|z^n\frac{\iota^n}{z^n}(1+\sum\limits_{j=1}^{\lfloor{\frac{n}{2}}\rfloor} (\iota)^{-2j} c_{2j}(S) z^{2j})|dz\\&= \frac{1}{\pi} \int\limits_{-\infty}^{\infty}\frac{1}{z^2} \log|\iota^n(1
 +\sum\limits_{j=1}^{\lfloor{\frac{n}{2}}\rfloor}(-1)^j c_{2j}(S) z^{2j})|dz
 \\&=\frac{1}{\pi} \int\limits_{-\infty}^{\infty}\frac{1}{z^2} \log|1+\sum\limits_{j=1}^{\lfloor{\frac{n}{2}}\rfloor} (-1)^jc_{2j}(S) z^{2j}|dz.
\end{align*}\qed

\noindent {\bf Remark 3.3} $(I)$. We note that the same integral expression holds for all non bipartite sidigraphs which have characteristic polynomial of the form $(2)$. It remains a problem to define a quasi-order relation (if possible) for sidigraphs in $\Delta^2_n$ for comparison of energy.\\
$(II)$. Let $D$ be a bipartite digraph on $n$ vertices and let $S_1$ and $S_2$ be sidigraphs on $D$ such that $S_1\in \Delta^1_n$ and $S_2\in \Delta^2_n$. Then by Theorems $2.5$ and $2.7$, the characteristic polynomials of $S_1$ and $S_2$ are given by\\
$\phi_{{S_1}(z)}=z^n+\sum_{j=1}^{\lfloor{\frac{n}{2}}\rfloor} (-1)^jc_{2j}(S_1)z^{n-2j}$ and $\phi_{{S_2}(z)}=z^n+\sum_{j=1}^{\lfloor{\frac{n}{2}}\rfloor} c_{2j}(S_2)z^{n-2j}$, where $c_{2j}(S_i)=|\pounds_{2j}(S_1)|=|\pounds_{2j}(S_2)|$ for $i=1,2$.\\
If $n$ is odd, then $\phi_{S_1}(\iota z)=\pm\iota\phi_{S_2}(z)$. The sign is $+$ or $-$ according as $n\equiv 1\pmod 4 $  or $n\equiv 3\pmod 4 $ respectively. If $n$ is even, then $\phi_{S_1}(\iota z)=\pm\phi_{S_2}(z)$. The sign is $+$ or $-$ according as $n\equiv 0\pmod 4 $ or $n\equiv 2\pmod 4 $ respectively.\\
Therefore, $spec(S_1)=\iota spec(S_2)$. From this relation energy of both sidigraphs can be calculated from the eigenvalues of the either sidigraph. For $i=1,2$, let $z_{ij}=\Re z_{ij}+\iota\Im z_{ij}$, where $j=1,2,\cdots,n$ be the eigenvalues of $S_i$ and $\Im z_j$ denotes the imaginary part of $z_j$. Then \\
$E(S_1)=\sum_{j=1}^n|\Re z_{1j}|=\sum_{j=1}^n|\Im z_{2j}|$ and $E(S_2)=\sum_{j=1}^n|\Re z_{2j}|=\sum_{j=1}^n|\Im z_{1j}|$.\\

In \cite{rgc}, the authors proved that the energy of a digraph in $\Delta_n$ decreases when we delete an arc from a cycle of length $2$. As in digraphs, in general it is not possible to predict the change in the energy of a non cycle balanced sidigraph by deleting an arc from a cycle of length $2$. It can decrease, increase or remain same by deleting an arc of a cycle of length $2$ as can be seen in the following example.\\

%TeXCAD Picture [f6.pic]. Options:
%\grade{\on}
%\emlines{\off}
%\epic{\off}
%\beziermacro{\on}
%\reduce{\on}
%\snapping{\off}
%\quality{8.00}
%\graddiff{0.01}
%\snapasp{1}
%\zoom{4.0000}
\unitlength 1mm % = 2.85pt
\linethickness{0.4pt}
\ifx\plotpoint\undefined\newsavebox{\plotpoint}\fi % GNUPLOT compatibility
\begin{picture}(111.75,50.38)(0,0)
\put(5.5,40.25){\circle*{2}}
\put(22.75,40){\circle*{2}}
\put(40,40){\circle*{2}}
\put(40,17){\circle*{2}}
\put(23.25,17){\circle*{2}}
\put(5.25,17){\circle*{2}}
\put(56.25,35.25){\circle*{2}}
\put(56.25,17.25){\circle*{2}}
\put(75.25,17.25){\circle*{2}}
\put(73,47.75){\circle*{2}}
\put(105.5,42){\circle*{2}}
\put(94.5,47.5){\circle*{2}}
\put(106.75,34.5){\circle*{2}}
\put(87.5,34){\circle*{2}}
\put(87.5,17.25){\circle*{2}}
\put(106.75,17){\circle*{2}}
\put(63.25,41.5){\circle*{2}}
%\vector[middle]{dot}(5.25,40.25)(22.75,40.25)
\put(14,40.25){\vector(1,0){.07}}\multiput(5.18,40.18)(.9722,0){19}{{\rule{.4pt}{.4pt}}}
%\end
%\vector[middle]{dot}(22.75,40.25)(39.75,40.25)
\put(31.25,40.25){\vector(1,0){.07}}\multiput(22.68,40.18)(.9444,0){19}{{\rule{.4pt}{.4pt}}}
%\end
%\vector[middle]{dot}(40,17)(40,40)
\put(40,28.5){\vector(0,1){.07}}\multiput(39.93,16.93)(0,.9583){25}{{\rule{.4pt}{.4pt}}}
%\end
%\vector[middle](5.25,40.5)(5.25,17.25)
\put(5.25,28.88){\vector(0,-1){.07}}\put(5.25,40.5){\line(0,-1){23.25}}
%\end
%\vector[middle](5.25,17.25)(23.25,17.25)
\put(14.25,17.25){\vector(1,0){.07}}\put(5.25,17.25){\line(1,0){18}}
%\end
%\vector[middle](23.25,17.25)(40.25,17.25)
\put(31.75,17.25){\vector(1,0){.07}}\put(23.25,17.25){\line(1,0){17}}
%\end
%\qbezvec[middle](39.75,40.25)(29.75,47.38)(22.75,40)
\put(30.5,43.75){\vector(-1,0){.07}}\qbezier(39.75,40.25)(29.75,47.38)(22.75,40)
%\end
%\qbezvec[middle](22.75,40.25)(15,47.13)(5.25,40.5)
\put(14.5,43.75){\vector(-1,0){.07}}\qbezier(22.75,40.25)(15,47.13)(5.25,40.5)
%\end
%\qbezvec[middle](72.75,48)(63.63,49.63)(63,41.75)
\put(65.75,47.25){\vector(-2,-1){.07}}\qbezier(72.75,48)(63.63,49.63)(63,41.75)
%\end
%\qbezvec[middle](105.25,42.5)(103.63,50.38)(94.5,47.75)
\put(101.75,47.75){\vector(-3,2){.07}}\qbezier(105.25,42.5)(103.63,50.38)(94.5,47.75)
%\end
%\vector[middle]{dot}(63,41.75)(72.75,48.25)
\put(67.88,45){\vector(3,2){.07}}\multiput(62.93,41.68)(.8125,.5417){13}{{\rule{.4pt}{.4pt}}}
%\end
%\vector[middle]{dot}(87.5,34.25)(106.75,34.75)
\put(97.13,34.5){\vector(1,0){.07}}\multiput(87.43,34.18)(.9625,.025){21}{{\rule{.4pt}{.4pt}}}
%\end
%\vector[middle](56,35.25)(62.5,41.75)
\put(59.25,38.5){\vector(1,1){.07}}\multiput(56,35.25)(.03367876,.03367876){193}{\line(0,1){.03367876}}
%\end
%\vector[middle](56,17.25)(56,35.5)
\put(56,26.38){\vector(0,1){.07}}\put(56,17.25){\line(0,1){18.25}}
%\end
%\vector[middle](56.5,35.75)(75.5,35.75)
\put(66,35.75){\vector(1,0){.07}}\put(56.5,35.75){\line(1,0){19}}
%\end
%\vector[middle](75.5,35.75)(75.5,17.75)
\put(75.5,26.75){\vector(0,-1){.07}}\put(75.5,35.75){\line(0,-1){18}}
%\end
%\vector[middle](75.25,17.25)(55.25,17.25)
\put(65.25,17.25){\vector(-1,0){.07}}\put(75.25,17.25){\line(-1,0){20}}
%\end
%\vector[middle](106.75,34.5)(106.75,17.25)
\put(106.75,25.88){\vector(0,-1){.07}}\put(106.75,34.5){\line(0,-1){17.25}}
%\end
%\vector[middle](106.5,17)(87.25,17)
\put(96.88,17){\vector(-1,0){.07}}\put(106.5,17){\line(-1,0){19.25}}
%\end
%\vector[middle](87.25,17.5)(87.25,34.5)
\put(87.25,26){\vector(0,1){.07}}\put(87.25,17.5){\line(0,1){17}}
%\end
%\vector[middle](94.25,47.5)(105.25,42)
\put(99.75,44.75){\vector(2,-1){.07}}\multiput(94.25,47.5)(.06707317,-.03353659){164}{\line(1,0){.06707317}}
%\end
%\vector[middle]{dot}(105.25,42.25)(106.5,34.75)
\put(105.88,38.5){\vector(1,-4){.07}}\multiput(105.18,42.18)(.1563,-.9375){9}{{\rule{.4pt}{.4pt}}}
%\end
\put(75.75,35.5){\circle*{1.8}}
\put(22,10.75){\makebox(0,0)[cc]{$S_1$}}
\put(65.75,11){\makebox(0,0)[cc]{$S_2$}}
\put(96.75,11){\makebox(0,0)[cc]{$S_3$}}
\put(56,6.25){\makebox(0,0)[cc]{Fig. $8$}}
\put(4,44.25){\makebox(0,0)[cc]{$v_1$}}
\put(22.75,35){\makebox(0,0)[cc]{$v_2$}}
\put(41,45.75){\makebox(0,0)[cc]{$v_3$}}
\put(39.75,11.75){\makebox(0,0)[cc]{$v_4$}}
\put(23,22.75){\makebox(0,0)[cc]{$v_5$}}
\put(5,11.75){\makebox(0,0)[cc]{$v_6$}}
\put(77.5,50){\makebox(0,0)[cc]{$v_1$}}
\put(60,44.75){\makebox(0,0)[cc]{$v_2$}}
\put(51.5,35.75){\makebox(0,0)[cc]{$v_3$}}
\put(51.25,17.25){\makebox(0,0)[cc]{$v_4$}}
\put(76,12.25){\makebox(0,0)[cc]{$v_5$}}
\put(75.75,39.75){\makebox(0,0)[cc]{$v_6$}}
\put(90.5,44.75){\makebox(0,0)[cc]{$v_1$}}
\put(110,44){\makebox(0,0)[cc]{$v_2$}}
\put(111.5,35){\makebox(0,0)[cc]{$v_3$}}
\put(111.75,16.75){\makebox(0,0)[cc]{$v_4$}}
\put(87,12.25){\makebox(0,0)[cc]{$v_5$}}
\put(84.75,38){\makebox(0,0)[cc]{$v_6$}}
\end{picture}

\noindent {\bf Example 3.4.} Consider the sidigraphs $S_1$, $S_2$ and $S_3$ as shown in Fig. $8$. It is easy to see that $\phi_{S_1}(z)=z^6+2z^4+1$ and $\phi_{S^{(v_1,v_2)}_1}(z)=z^6+z^4+1$, where $S_1^{(v_1,v_2)}$ denotes the sidigraph obtained by deleting the arc $(v_1,v_2)$. Note $E(S_1)\approx 2.4916$ and $E(S_1^{(v_1,v_2)})\approx 2.9104$. So the energy increases in this case. Also, $\phi_ {S_2}(z)=z^6+z^4-z^2-1$ and $spec(S_2)=\{-1,1,-\iota^{(2)},\iota^{(2)}\}$ so that $E(S_2)=2$. If we delete arc $(v_1,v_2)$, the resulting sidigraph has eigenvalues $\{-1,0^{(2)},1,-\iota, \iota\}$ so the energy of the resulting sidigraph is again $2$. That is, the energy remains same in this case. It is not difficult to check that $E(S_3)=2+2\sqrt{2}$ and $E(S_3^{(v_1,v_2)})=2\sqrt{2}$. So the energy decreases in this case.\\

The following result shows that the energy of a sidigraph in $\Delta^1_n$ decreases when we delete an arc from a cycle of length $2$. The proof is same as the proof of [Theorem $2.6$, \cite{rgc}].\\

\noindent {\bf Theorem 3.5.} Let $S$ be a sidigraph in $\Delta^1_n$ with a pair of symmetric arcs and let $S^\prime$ be the sidigraph obtained by deleting one of these arcs. Then $E(S^\prime)<E(S)$.\\

\section{Equienergetic sidigraphs}

Two noncospectral sidigraphs $S_1$ and $S_2$ are said to be equienergetic if $E(S_1)=E(S_2)$. For equienergetic (si)digraphs see \cite{lr,pm}. We construct a new family of strongly connected equienergetic sidigraphs and answer to the open problem $2$ raised in \cite{pm}.\\
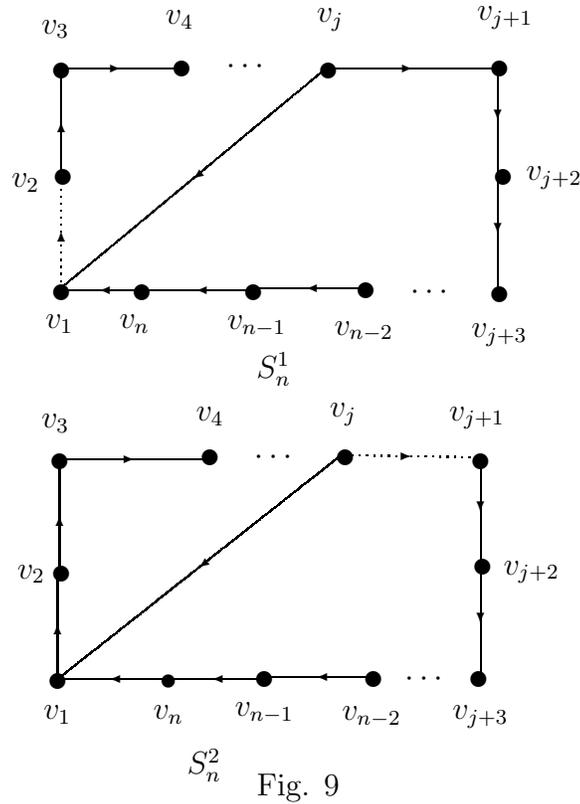
\begin{align*}
%TeXCAD Picture [f7.pic]. Options:
%\grade{\on}
%\emlines{\off}
%\epic{\off}
%\beziermacro{\on}
%\reduce{\on}
%\snapping{\off}
%\quality{8.00}
%\graddiff{0.01}
%\snapasp{1}
%\zoom{4.0000}
\unitlength 1mm % = 2.85pt
\linethickness{0.4pt}
\ifx\plotpoint\undefined\newsavebox{\plotpoint}\fi % GNUPLOT compatibility
\begin{picture}(74.25,104.75)(0,0)
\put(8.25,16.5){\circle*{2}}
\put(64.25,16.75){\circle*{2}}
\put(8.75,30.75){\circle*{2}}
\put(8.5,45.75){\circle*{2}}
\put(64.5,45.75){\circle*{2}}
\put(64.75,31.75){\circle*{2}}
%\vector[middle](64.5,46)(64.5,32)
\put(64.5,39){\vector(0,-1){.07}}\put(64.5,46){\line(0,-1){14}}
%\end
%\vector[middle](64.5,32)(64.5,16.5)
\put(64.5,24.25){\vector(0,-1){.07}}\put(64.5,32){\line(0,-1){15.5}}
%\end
\put(28.5,46.25){\circle*{2}}
\put(46.5,46.25){\circle*{2}}
\put(50.25,16.75){\circle*{2}}
\put(35.75,16.75){\circle*{2}}
%\vector[middle](8.5,31)(8.5,46)
\put(8.5,38.5){\vector(0,1){.07}}\put(8.5,31){\line(0,1){15}}
%\end
%\vector[middle](8.5,46)(28.25,46)
\put(18.38,46){\vector(1,0){.07}}\put(8.5,46){\line(1,0){19.75}}
%\end
%\vector[middle](50,17)(35.5,17)
\put(42.75,17){\vector(-1,0){.07}}\put(50,17){\line(-1,0){14.5}}
%\end
%\vector[middle](35.25,16.75)(22.5,16.75)
\put(28.88,16.75){\vector(-1,0){.07}}\put(35.25,16.75){\line(-1,0){12.75}}
%\end
%\vector[middle](22.5,16.75)(8,16.75)
\put(15.25,16.75){\vector(-1,0){.07}}\put(22.5,16.75){\line(-1,0){14.5}}
%\end
\put(23,16.5){\circle*{1.8}}
\put(8.25,12){\makebox(0,0)[cc]{$v_1$}}
\put(4.75,30.75){\makebox(0,0)[cc]{$v_2$}}
\put(8,50.75){\makebox(0,0)[cc]{$v_3$}}
\put(28.5,51.75){\makebox(0,0)[cc]{$v_4$}}
\put(46.25,52){\makebox(0,0)[cc]{$v_j$}}
\put(64,51.5){\makebox(0,0)[cc]{$v_{j+1}$}}
\put(64.5,11.75){\makebox(0,0)[cc]{$v_{j+3}$}}
\put(50,11.75){\makebox(0,0)[cc]{$v_{n-2}$}}
\put(35.75,12.25){\makebox(0,0)[cc]{$v_{n-1}$}}
\put(23,12){\makebox(0,0)[cc]{$v_{n}$}}
\put(37.25,46.25){\makebox(0,0)[cc]{$\cdots$}}
\put(57.25,16.75){\makebox(0,0)[cc]{$\cdots$}}
\put(71.25,31){\makebox(0,0)[cc]{$v_{j+2}$}}
%\vector[middle](46.25,46.75)(8,16.5)
\put(27.13,31.63){\vector(-4,-3){.07}}\multiput(46.25,46.75)(-.0426421405,-.0337235229){897}{\line(-1,0){.0426421405}}
%\end
%\vector[middle](8.25,17.25)(8.25,30.75)
\put(8.25,24){\vector(0,1){.07}}\put(8.25,17.25){\line(0,1){13.5}}
%\end
%\vector[middle]{dot}(46.5,46.5)(64,46.25)
\put(55.25,46.38){\vector(1,0){.07}}\multiput(46.43,46.43)(.9722,-.0139){19}{{\rule{.4pt}{.4pt}}}
%\end
\put(8.75,68.25){\circle*{2}}
\put(9,83.5){\circle*{2}}
\put(8.75,97.75){\circle*{2}}
\put(24.75,98){\circle*{2}}
\put(44.25,97.75){\circle*{2}}
\put(67,98){\circle*{2}}
\put(67,68){\circle*{2}}
\put(67.5,83.5){\circle*{2}}
\put(49.25,68.5){\circle*{2}}
\put(34.25,68.25){\circle*{2}}
\put(19.5,68.25){\circle*{2}}
%\vector[middle](8.75,83.75)(8.75,98)
\put(8.75,90.88){\vector(0,1){.07}}\put(8.75,83.75){\line(0,1){14.25}}
%\end
%\vector[middle](8.5,98)(24.75,98)
\put(16.63,98){\vector(1,0){.07}}\put(8.5,98){\line(1,0){16.25}}
%\end
%\vector[middle](44,98)(66.75,98)
\put(55.38,98){\vector(1,0){.07}}\put(44,98){\line(1,0){22.75}}
%\end
%\vector[middle](66.75,98.25)(66.75,84)
\put(66.75,91.13){\vector(0,-1){.07}}\put(66.75,98.25){\line(0,-1){14.25}}
%\end
%\vector[middle](66.75,83.75)(66.75,68.25)
\put(66.75,76){\vector(0,-1){.07}}\put(66.75,83.75){\line(0,-1){15.5}}
%\end
%\vector[middle](49,68.75)(33.75,68.75)
\put(41.38,68.75){\vector(-1,0){.07}}\put(49,68.75){\line(-1,0){15.25}}
%\end
%\vector[middle](34.25,68.5)(19.25,68.5)
\put(26.75,68.5){\vector(-1,0){.07}}\put(34.25,68.5){\line(-1,0){15}}
%\end
%\vector[middle](19.25,68.5)(8.5,68.5)
\put(13.88,68.5){\vector(-1,0){.07}}\put(19.25,68.5){\line(-1,0){10.75}}
%\end
%\vector[middle]{dot}(8.75,68.75)(8.75,83.75)
\put(8.75,76.25){\vector(0,1){.07}}\multiput(8.68,68.68)(0,.9375){17}{{\rule{.4pt}{.4pt}}}
%\end
%\vector[middle](44,98)(8.75,68.75)
\put(26.38,83.38){\vector(-4,-3){.07}}\multiput(44,98)(-.0406574394,-.0337370242){867}{\line(-1,0){.0406574394}}
%\end
\put(33.5,98.25){\makebox(0,0)[cc]{$\cdots$}}
\put(58,68.25){\makebox(0,0)[cc]{$\cdots$}}
\put(40.25,2.25){\makebox(0,0)[cc]{Fig. $9$}}
\put(8.5,64){\makebox(0,0)[cc]{$v_1$}}
\put(4,82.75){\makebox(0,0)[cc]{$v_2$}}
\put(8,103.25){\makebox(0,0)[cc]{$v_3$}}
\put(24.5,104.5){\makebox(0,0)[cc]{$v_4$}}
\put(44.5,104.5){\makebox(0,0)[cc]{$v_j$}}
\put(67.75,104.75){\makebox(0,0)[cc]{$v_{j+1}$}}
\put(74.25,83.5){\makebox(0,0)[cc]{$v_{j+2}$}}
\put(67,62.5){\makebox(0,0)[cc]{$v_{j+3}$}}
\put(49,63){\makebox(0,0)[cc]{$v_{n-2}$}}
\put(34.75,63.75){\makebox(0,0)[cc]{$v_{n-1}$}}
\put(18.5,64){\makebox(0,0)[cc]{$v_n$}}
\put(37,58){\makebox(0,0)[cc]{$S^1_n$}}
\put(27.75,5.25){\makebox(0,0)[cc]{$S^2_n$}}
\end{picture}
\end{align*}

\noindent {\bf Theorem 4.1.} For each positive integer $n\ge 4$, there exists a pair of  noncospectral, equienergetic and strongly connected sidigraphs of order $n$ with both the constituents non cycle balanced.\\
\noindent {\bf Proof.} {\bf Case $1$.} $n$ is even. For each even positive integer $n\ge 4$, consider the sidigraphs $S^1_n$ and $S^2_n$ with vertex and arc sets given by \\
$V(S^1_n)=V(S^2_n)=\{v_{1},v_{2},\cdots,v_{n}\}$,\\ $\mathscr{A}(S^1_n)=\{[v_1,v_2],(v_2,v_3),\cdots,(v_{j-1},v_j),(v_j,v_{j+1}),\cdots,(v_n,v_1),(v_j,v_1)\}$\\
and \\
$\mathscr{A}(S^2_n)=\{(v_1,v_2),\cdots,(v_{j-1},v_j),[v_j,v_{j+1}],(v_{j+1},v_{j+2})\cdots,(v_n,v_1),(v_j,v_1)\}$,\\
where symbol $[u,v]$ means sign of arc from vertex $u$ to vertex $v$ is negative and we choose vertex $v_j$ such that positive integer $j$ is odd. The sidigraphs so constructed are shown in Fig. $9$. Sidigraph $S^1_n$ has one negative cyclic arc $[v_1,v_2]$ and sidigraph $S^2_n$ also has one negative cyclic arc $[v_j,v_{j+1}]$. Therefore both sidigraphs are non cycle balanced for each even $n$.\\
\indent By Theorem $1.1$, $\phi_{S^1_n}(z)=z^n+z^{n-j}+1$ and $\phi_{S^2_n}(z)=z^n-z^{n-j}+1$.\\
Clearly $S^1_n$ and $S^2_n$ are noncospectral. As $n$ is even and $j$ is odd, so $n-j$ is odd.\\
Therefore, $\phi_{S^1_n}(-z)=z^n-z^{n-j}+1=\phi_{S^2_n}(z)$ and hence $spec(S^1_n)=-spec(S^2_n)$.\\
Thus, $E(S^1_n)=E(S^2_n).$ \\

\begin{align*}
%TeXCAD Picture [f8.pic]. Options:
%\grade{\on}
%\emlines{\off}
%\epic{\off}
%\beziermacro{\on}
%\reduce{\on}
%\snapping{\off}
%\quality{8.00}
%\graddiff{0.01}
%\snapasp{1}
%\zoom{4.0000}
\unitlength 1mm % = 2.85pt
\linethickness{0.4pt}
\ifx\plotpoint\undefined\newsavebox{\plotpoint}\fi % GNUPLOT compatibility
\begin{picture}(74.25,104.75)(0,0)
\put(8.25,16.5){\circle*{2}}
\put(64.25,16.75){\circle*{2}}
\put(8.75,30.75){\circle*{2}}
\put(8.5,45.75){\circle*{2}}
\put(64.5,45.75){\circle*{2}}
\put(64.75,31.75){\circle*{2}}
%\vector[middle](64.5,46)(64.5,32)
\put(64.5,39){\vector(0,-1){.07}}\put(64.5,46){\line(0,-1){14}}
%\end
%\vector[middle](64.5,32)(64.5,16.5)
\put(64.5,24.25){\vector(0,-1){.07}}\put(64.5,32){\line(0,-1){15.5}}
%\end
\put(28.5,46.25){\circle*{2}}
\put(46.5,46.25){\circle*{2}}
\put(50.25,16.75){\circle*{2}}
\put(35.75,16.75){\circle*{2}}
%\vector[middle](8.5,31)(8.5,46)
\put(8.5,38.5){\vector(0,1){.07}}\put(8.5,31){\line(0,1){15}}
%\end
%\vector[middle](8.5,46)(28.25,46)
\put(18.38,46){\vector(1,0){.07}}\put(8.5,46){\line(1,0){19.75}}
%\end
%\vector[middle](50,17)(35.5,17)
\put(42.75,17){\vector(-1,0){.07}}\put(50,17){\line(-1,0){14.5}}
%\end
%\vector[middle](35.25,16.75)(22.5,16.75)
\put(28.88,16.75){\vector(-1,0){.07}}\put(35.25,16.75){\line(-1,0){12.75}}
%\end
%\vector[middle](22.5,16.75)(8,16.75)
\put(15.25,16.75){\vector(-1,0){.07}}\put(22.5,16.75){\line(-1,0){14.5}}
%\end
\put(23,16.5){\circle*{1.8}}
\put(8.25,12){\makebox(0,0)[cc]{$v_1$}}
\put(4.75,30.75){\makebox(0,0)[cc]{$v_2$}}
\put(8,50.75){\makebox(0,0)[cc]{$v_3$}}
\put(28.5,51.75){\makebox(0,0)[cc]{$v_4$}}
\put(46.25,52){\makebox(0,0)[cc]{$v_j$}}
\put(64,51.5){\makebox(0,0)[cc]{$v_{j+1}$}}
\put(64.5,11.75){\makebox(0,0)[cc]{$v_{j+3}$}}
\put(50,11.75){\makebox(0,0)[cc]{$v_{n-3}$}}
\put(35.75,12.25){\makebox(0,0)[cc]{$v_{n-2}$}}
\put(23,12){\makebox(0,0)[cc]{$v_{n-1}$}}
\put(37.25,46.25){\makebox(0,0)[cc]{$\cdots$}}
\put(57.25,16.75){\makebox(0,0)[cc]{$\cdots$}}
\put(71.25,31){\makebox(0,0)[cc]{$v_{j+2}$}}
%\vector[middle](46.25,46.75)(8,16.5)
\put(27.13,31.63){\vector(-4,-3){.07}}\multiput(46.25,46.75)(-.0426421405,-.0337235229){897}{\line(-1,0){.0426421405}}
%\end
%\vector[middle]{dot}(46.5,46.5)(64,46.25)
\put(55.25,46.38){\vector(1,0){.07}}\multiput(46.43,46.43)(.9722,-.0139){19}{{\rule{.4pt}{.4pt}}}
%\end
\put(8.75,68.25){\circle*{2}}
\put(9,83.5){\circle*{2}}
\put(8.75,97.75){\circle*{2}}
\put(24.75,98){\circle*{2}}
\put(44.25,97.75){\circle*{2}}
\put(67,98){\circle*{2}}
\put(67,68){\circle*{2}}
\put(49.25,68.5){\circle*{2}}
\put(34.25,68.25){\circle*{2}}
\put(19.5,68.25){\circle*{2}}
%\vector[middle](8.5,98)(24.75,98)
\put(16.63,98){\vector(1,0){.07}}\put(8.5,98){\line(1,0){16.25}}
%\end
%\vector[middle](44,98)(66.75,98)
\put(55.38,98){\vector(1,0){.07}}\put(44,98){\line(1,0){22.75}}
%\end
%\vector[middle](66.75,98.25)(66.75,84)
\put(66.75,91.13){\vector(0,-1){.07}}\put(66.75,98.25){\line(0,-1){14.25}}
%\end
%\vector[middle](66.75,83.75)(66.75,68.25)
\put(66.75,76){\vector(0,-1){.07}}\put(66.75,83.75){\line(0,-1){15.5}}
%\end
%\vector[middle](49,68.75)(33.75,68.75)
\put(41.38,68.75){\vector(-1,0){.07}}\put(49,68.75){\line(-1,0){15.25}}
%\end
%\vector[middle](34.25,68.5)(19.25,68.5)
\put(26.75,68.5){\vector(-1,0){.07}}\put(34.25,68.5){\line(-1,0){15}}
%\end
%\vector[middle](19.25,68.5)(8.5,68.5)
\put(13.88,68.5){\vector(-1,0){.07}}\put(19.25,68.5){\line(-1,0){10.75}}
%\end
%\vector[middle]{dot}(8.75,68.75)(8.75,83.75)
\put(8.75,76.25){\vector(0,1){.07}}\multiput(8.68,68.68)(0,.9375){17}{{\rule{.4pt}{.4pt}}}
%\end
%\vector[middle](44,98)(8.75,68.75)
\put(26.38,83.38){\vector(-4,-3){.07}}\multiput(44,98)(-.0406574394,-.0337370242){867}{\line(-1,0){.0406574394}}
%\end
\put(33.5,98.25){\makebox(0,0)[cc]{$\cdots$}}
\put(58,68.25){\makebox(0,0)[cc]{$\cdots$}}
\put(40.25,2.25){\makebox(0,0)[cc]{Fig. $10$}}
\put(8.5,64){\makebox(0,0)[cc]{$v_1$}}
\put(4,82.75){\makebox(0,0)[cc]{$v_2$}}
\put(8,103.25){\makebox(0,0)[cc]{$v_3$}}
\put(24.5,104.5){\makebox(0,0)[cc]{$v_4$}}
\put(44.5,104.5){\makebox(0,0)[cc]{$v_j$}}
\put(67.75,104.75){\makebox(0,0)[cc]{$v_{j+1}$}}
\put(74.25,83.5){\makebox(0,0)[cc]{$v_{j+2}$}}
\put(67,62.5){\makebox(0,0)[cc]{$v_{j+3}$}}
\put(37,58){\makebox(0,0)[cc]{$S^3_n$}}
\put(27.75,5.25){\makebox(0,0)[cc]{$S^4_n$}}
\put(66.75,83.75){\circle*{1.8}}
\put(49.5,63.25){\makebox(0,0)[cc]{$v_{n-3}$}}
\put(34.25,64.25){\makebox(0,0)[cc]{$v_{n-2}$}}
\put(19.25,63.75){\makebox(0,0)[cc]{$v_{n-1}$}}
%\vector[middle](8.75,31)(22,35.5)
\put(15.38,33.25){\vector(3,1){.07}}\multiput(8.75,31)(.0988806,.03358209){134}{\line(1,0){.0988806}}
%\end
%\vector[middle](21.5,35.25)(8,16.75)
\put(14.75,26){\vector(-3,-4){.07}}\multiput(21.5,35.25)(-.033665835,-.046134663){401}{\line(0,-1){.046134663}}
%\end
%\vector[middle](8.75,84)(22.5,89.75)
\put(15.63,86.88){\vector(2,1){.07}}\multiput(8.75,84)(.08040936,.03362573){171}{\line(1,0){.08040936}}
%\end
%\vector[middle](22.25,89.25)(8.75,69)
\put(15.5,79.13){\vector(-2,-3){.07}}\multiput(22.25,89.25)(-.033665835,-.050498753){401}{\line(0,-1){.050498753}}
%\end
\put(22.5,89.75){\circle*{2}}
\put(21.75,35.25){\circle*{2}}
%\vector[middle](8.25,16.75)(8.25,31)
\put(8.25,23.88){\vector(0,1){.07}}\put(8.25,16.75){\line(0,1){14.25}}
%\end
%\vector[middle](8.75,83.5)(8.75,98)
\put(8.75,90.75){\vector(0,1){.07}}\put(8.75,83.5){\line(0,1){14.5}}
%\end
\put(26.75,91){\makebox(0,0)[cc]{$v_n$}}
\put(26.25,37.25){\makebox(0,0)[cc]{$v_n$}}
\end{picture}
\end{align*}

{\bf Case $2$.} $n$ is odd.  For each odd positive integer $n\ge 5$, consider the sidigraphs $S^3_n$ and $S^4_n$ with vertex and arc sets given by \\
$V(S^3_n)=V(S^4_n)=\{v_{1},v_{2},\cdots,v_{n}\}$,\\ $\mathscr{A}(S^3_n)=\{[v_1,v_2],(v_2,v_3),\cdots,(v_{j-1},v_j),(v_j,v_{j+1}),\cdots,(v_{n-1},v_1),(v_j,v_1),\\(v_2,v_n),(v_n,v_1)\}$\\
and \\
$\mathscr{A}(S^4_n)=\{(v_1,v_2),\cdots,(v_{j-1},v_j),[v_j,v_{j+1}],(v_{j+1},v_{j+2})\cdots,(v_{n-1},v_1),(v_j,v_1),\\(v_2,v_n),(v_n,v_1)\}$.\\

The sidigraphs so constructed are shown in Fig. $10$, where  for vertex $v_j$, positive integer $j$ is odd. Clearly, both the sidigraphs are non cycle balanced. By Theorem $1.1$, we have $\phi_{S^3_n}(z)=z^n+z^{n-3}+z^{n-j}+z$ and $\phi_{S^4_n}(z)=z^n-z^{n-3}-z^{n-j}+z$. Clearly, $S^3_n$ and $S^4_n$ are noncospectral. As both $n$ and $j$ are odd, so both $n-3$ and $n-j$ are even. Also, $\phi_{S^3_n}(-z)=-z^n+z^{n-3}+z^{n-j}-z=-\phi_{S^4_n}(z)$. Therefore, $spec(S^3_n)=-spec(S^4_n)$ and so $E(S^3_n)=E(S^4_n)$.\qed

\section{Open problem}
We conclude this paper with the following open problem for further study.\\

{\bf Problem 5.1.} Characterize sidigraphs with characteristic polynomial of the form $(2)$ or $(3)$.\\

{\bf Acknowledgements}. The first author thanks University Grants Commission, New Delhi, India for providing senior research fellowship.\\

\end{document}